\documentclass[12pt,reqno]{amsart}
\usepackage{amssymb}
\usepackage{graphicx}
\usepackage{amsmath}
\usepackage{a4wide}
\usepackage{epsf}
\usepackage{psfrag}
\newtheorem{theorem}{Theorem}

\newtheorem{corollary}[theorem]{Corollary}
\newtheorem{lemma}[theorem]{Lemma}
\newtheorem{proposition}[theorem]{Proposition}

\theoremstyle{remark}

\def\gra{\operatorname{graph}}
\def\leq{\leqslant}
\def\geq{\geqslant}
\def\al{\alpha}
\def\be{\beta}

\def\la{\lambda}

\def\N{\mathbb N}
\def\R{\mathbb R}

\def\Z{\mathbb Z}

\def\cc {,\dots,}

\def\pq{q^{\phantom{q}}}      % to make subscripts lower!
\def\pchi{\chi^{\phantom{q}}} % to make subscripts lower!

\def\bei{\begin{itemize}}
\def\bee{\begin{enumerate}}
\def\ene{\end{enumerate}}
\def\eni{\end{itemize}}

\begin{document}
\title[Salem, Pisot, Mahler and graphs]
{Salem numbers, Pisot numbers, Mahler measure and graphs}
\author{ James McKee}
\address{Department of Mathematics\\
Royal Holloway, University of London\\
Egham Hill\\
Egham\\
Surrey TW20 0EX\\
UK}
\email{James.McKee@rhul.ac.uk}
\author{ Chris Smyth}
\address{School of Mathematics \\
University of Edinburgh\\
James Clerk Maxwell Building\\
King's Buildings, Mayfield Road\\
Edinburgh EH9 3JZ\\
Scotland, U.K.} \email{C.Smyth@ed.ac.uk} 
\date{\today}

\maketitle

\begin{abstract} We use graphs to define  sets of Salem and Pisot numbers, and prove that the union of these sets is closed, supporting a conjecture of Boyd that the set of all Salem and Pisot numbers is closed. We find all trees that define Salem numbers. We show that for all integers $n$ the smallest known element of the $n$-th derived set of the set of Pisot numbers comes from
a graph. We define the Mahler measure of a graph, and find all graphs of Mahler measure less than $\frac12(1+\sqrt{5})$. Finally, we list all small Salem numbers known to be definable using a graph.
\end{abstract}

\section{Introduction}

The work described in this paper arose from the following idea:
that one way of studying algebraic integers might be by associating combinatorial
objects with them. Here, we try to do this for two particular classes of algebraic integers,
 Salem numbers and Pisot numbers, the associated combinatorial objects being graphs.
We also  find all graphs of small Mahler measure.
 All but one of these measures turns out to be a Salem number.

A {\em Pisot number} is a real algebraic integer $\theta>1$,
all of whose other Galois conjugates have modulus strictly
less than 1.
A {\em Salem number} is a real algebraic integer $\tau>1$,
whose other conjugates all have modulus at most 1, with at least
one having modulus exactly 1.
It follows that the minimal polynomial $P(z)$ of $\tau$ is
{\em reciprocal} (that is, $z^{{\rm deg\ }P}P(1/z)=P(z)$), that
$\tau^{-1}$ is a conjugate of $\tau$, that all conjugates of
$\tau$ other than $\tau$ and $\tau^{-1}$ have modulus exactly 1,
and that $P(z)$ has even degree.
The set of all Pisot numbers is traditionally (if a little
unfortunately) denoted $S$, with $T$ being used for the set of all
Salem numbers.

 We call a graph $G$ a {\it Salem graph} if either
\begin{itemize}
\item it is nonbipartite, has only one eigenvalue $\la>2$
  and no eigenvalues
in $(-\infty,-2)$;

or

\item it is bipartite, has only one eigenvalue $\la>2$ and only the
eigenvalue $-\la$ in $(-\infty,-2)$. \eni

We call a Salem graph {\it trivial} if it is nonbipartite and
$\la\in\Z$, or it is bipartite and $\la^2\in\Z$. For a
nontrivial Salem graph, its associated Salem number $\tau(G)$ is
then the larger root of $z+1/z=\la$ in the nonbipartite case, and of
$\sqrt{z}+1/\sqrt{z}=\la$ in the bipartite case. (Proposition
\ref{P-Salem} shows that $\tau(G)$ is indeed a Salem number.) We
call $\tau(G)$ a {\it graph Salem number}, and denote by $T_{\gra}$
 the set of all graph Salem numbers. (For a trivial Salem graph $G$,
 $\tau(G)$ is a reciprocal quadratic Pisot number.)

Our first result is the following.

\begin{theorem} \label{T-closed}
  The set of limit points of $T_{\gra}$
is some set $S_{\gra}$ of Pisot numbers. Furthermore,
$T_{\gra}\cup S_{\gra}$ is closed.
\end{theorem}

In \cite[Corollary 9]{MRS}, a construction was given for certain
subsets $S^*$ of $S$ and $T^*$ of $T$, using a restricted class of
graphs (star-like trees). We showed that
 $T^*$ had its limit points in  $S^*$, and that (like $S$) $S^*$ was closed in $\R$.
  The main aim of
this paper is to push these ideas as far as we can.

We call elements of $S_{\gra}$ {\em graph Pisot numbers}.
The proof of Theorem \ref{T-closed} reveals a way to represent graph Pisot
numbers by bi-vertex-coloured graphs, which we call {\em Pisot
graphs}.

Since Boyd has long conjectured that $S$ is the set of limit
points of $T$, and that therefore $S\cup T$ is closed (\cite{Bo}),
our result is a step in the direction of a proof of his
conjecture. However, we can find elements in $T-T_{\gra}$
(see Section \ref{S-SmSa}) and elements in $S-S_{\gra}$ (see
Corollary \ref{C-small-Pisot}), so that graphs do not tell the whole
story.

 It is clearly desirable to describe all Salem graphs. While we
 have not been able to do this completely, we are able in
 Proposition \ref{P-maxdeg3} to restrict
 the class of graphs that can be Salem graphs. Naturally enough, we call a Salem graph that
happens to be a tree a {\em Salem tree}.
 In Section \ref{S-trees} we
 completely describe all Salem trees.

In Section \ref{S-Bertin} we show that the smallest known elements
of the $k$-th derived set of $S$ belong to the $k$-th derived set of $S_{\gra}$.
In Section \ref{S-Mahler}, we find all graphs having Mahler measure
at most $\frac12(1+\sqrt{5})$.
Finally, in Section \ref{S-SmSa} we list some small Salem numbers coming from graphs.

{\it Acknowledgments.} We are very grateful to Peter Rowlinson for providing us
with many references on graph eigenvalues. We also thank the referees for
helpful comments.

\section{Lemmas on graph eigenvalues}

For a graph $G$, recall that its eigenvalues are defined to be
those of its adjacency matrix $A=(a_{ij})$, where $a_{ij}=1$ if the
$i$th and $j$th vertices are joined by an edge (`adjacent'), and $0$
otherwise. Because $A$ is symmetric, all eigenvalues of $G$ are
real.

 The following
facts are essential ingredients in our proofs.

\begin{lemma}[{Interlacing Theorem.
 See
\cite[Theorem 9.1.1]{GR}}]\label{L-interlacing}

If a graph $G$ has eigenvalues
$\lambda_1 \leq \ldots \leq \lambda_n$, and a vertex of $G$ is
deleted to produce a graph $H$ with eigenvalues $\mu_1 \leq \ldots
\leq \mu_{n-1}$, then the eigenvalues of $G$ and $H$ interlace,
namely
$$
\lambda_1\leq\mu_1\leq\lambda_2\leq\mu_2\leq\cdots\leq\mu_{n-1}\leq
\lambda_n\,.
$$
\end{lemma}

    We denote the largest eigenvalue of a graph $G$, called its {\em
index},  by $\lambda(G)$. We call a graph that has all its
eigenvalues in the interval $[-2,2\,]$ a {\em cyclotomic graph}.
Connected graphs that have index at most $2$ have been classified,
and in fact all are cyclotomic.

\begin{lemma} [{ \cite{Smi},
\cite{Neu}---see also \cite[Theorem 2.1]{CR} }] \label{L-except}
The connected cyclotomic graphs are precisely the induced
subgraphs of the   graphs $\tilde E_6, \tilde E_7$
and $ \tilde E_8$, and those of the $(n+1)$-vertex graphs $\tilde A_n$
$(n\geq 2)$,  $\tilde
D_n$  $(n\geq 4)$, as in Figure \ref{F-except}.
\end{lemma}

\begin{figure}[h]
\begin{center}
\leavevmode
\psfragscanon
\psfrag{A}[l]{$\tilde A_n$}
\psfrag{D}[l]{$\tilde D_n$}
\psfrag{6}[l]{$\tilde E_6$}
\psfrag{7}{$\tilde E_7$}
\psfrag{8}{$\tilde E_8$}
\includegraphics[scale=0.5]{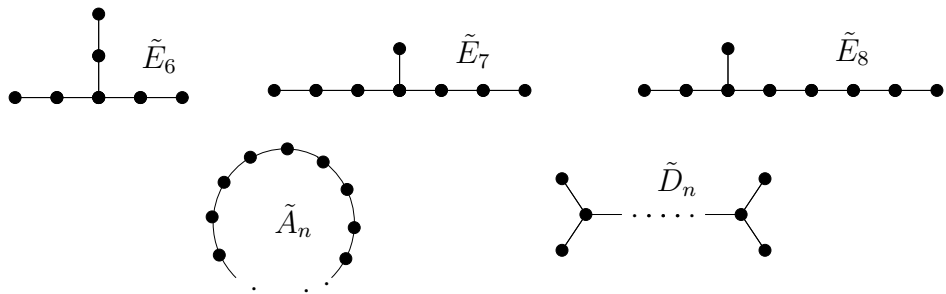}
\end{center}
\caption{The maximal connected cyclotomic graphs $\tilde E_6, \tilde E_7,
 \tilde E_8, \tilde A_n (n\geq 2)$ and $ \tilde
D_n (n\geq 4)$. The number of vertices is one more than the
subscript.}\label{F-except}
\end{figure}

Clearly, general cyclotomic graphs are then graphs all of whose
connected components are cyclotomic.

 An {\it internal path} of a graph $G$ is a sequence of vertices $x_1\cc
x_k$ of $G$ such that all vertices (except possibly $x_1$ and $x_k$) are
distinct, $x_i$ is adjacent to $x_{i+1}$ for $i=1\cc k-1$, $x_1$ and $x_k$ have
degree at least $3$, while $x_2\cc x_{k-1}$ have degree $2$. An
{\it internal edge} is an edge on an internal path.

\begin{lemma} \label{L-addedge}\begin{itemize}

\item[(i)] Suppose that the connected graph $G$ has $G'$ as a proper subgraph.
Then $\la(G')< \la(G)$.

\item[(ii)] Suppose that $G^*$ is a graph obtained from a connected graph
$G$ by subdividing an internal edge. Then $\la(G^*)\leq \la(G)$,
with equality if and only if $G=\tilde D_n$ for some $n\geq 5$.
\end{itemize}
\end{lemma}

For the proof, see \cite[Theorem 8.8.1(b)]{GR} and \cite[Proposition
2.4]{HS}. Note that on subdividing a (noninternal) edge of  $\tilde
A_n$,  $\la(G)=2$ does not change. For any other connected graph
$G$, if we subdivide a noninternal edge of $G$ to get a graph $G^*$,
then $G$ is (isomorphic to) a subgraph of $G^*$, so that, by (i),
$\la(G^*)> \la(G)$.

\begin{lemma}[{See \cite[Theorem 1.3]{CR} and references therein}]\label{L-maxdeg}
Every vertex of a graph $G$ has degree at most $\la(G)^2$.
\end{lemma}
\begin{proof} Suppose $G$ has a vertex of degree $d\geq 1$
(the result is trivial if every vertex has degree 0) Then the star
subgraph $G'$ of $G$ on that vertex and its adjacent vertices has
$\la(G')=\sqrt{d}$. This follows from the fact that its ``quotient''
(see Section \ref{S-quotients}) is $(z+1-zd/(z+1))^{-1}$, so that it
has two distinct eigenvalues $\pm\la$ satisfying
 $\la^2=(\sqrt{z}+1/\sqrt{z})^2=d$.
(If $d\geq 2$, then $0$ is also an eigenvalue.)
 By Lemma
\ref{L-addedge}(i), $\la(G')\leq \la(G)$,  giving the result.
\end{proof}

\section{Salem  graphs}

Let $G$ be a graph
on $n$ vertices, and let $\pchi_G(x)$ be its characteristic
polynomial (the characteristic polynomial of the adjacency matrix
of $G$).

When $G$ is nonbipartite, we define the {\em reciprocal
polynomial} of $G$, denoted $R_G(x)$,  by
$$
R_G(z)=z^{n}\pchi_G(z+1/z)\,.
$$
By construction $R_G$ is indeed a reciprocal polynomial, its roots
coming in pairs, each root $\beta=\al+1/\al$ of $\pchi_G$
corresponding to the (multiset) pair $\{\al,1/\al\}$ of roots of
$R_G$.

 When  $G$ is bipartite, the reciprocal polynomial $R_G(x)$ is defined by
$$
R_G(z)=z^{n/2}\pchi_G(\sqrt{z}+1/\sqrt{z})\,.
$$
In this case
$\pchi_G(-x)=(-1)^n\pchi_G(x)$: the characteristic polynomial is
either even or odd. From this one readily sees that $R_G(z)$ is
indeed a polynomial, and the correspondence this time is between
 the pairs $\{\beta,-\beta\}$ and $\{\alpha,1/\alpha\}$, where $\beta=\sqrt{\alpha}+1/\sqrt{\alpha}$. (We may
suppose that the branch of the square root is chosen such that
$\beta\geq 0$).

 As the roots of
$\pchi_G$ are all real,  in both cases the roots of $R_G$ are either real, or lie on the unit circle: if $\beta>2$ then the above correspondence is with a pair $\{\alpha,1/\alpha\}$, both positive; if $\beta\in[-2,2\,]$ then it is with
the pair $\{\alpha,1/\alpha=\bar\al\}$, both of modulus $1$.

 \begin{proposition}\label{P-Salem}
For a cyclotomic graph $G$, $R_G$ is indeed a cyclotomic polynomial.
For a nontrivial Salem graph $G$,
 $\tau(G)$ is indeed a Salem number.
\end{proposition}
 \begin{proof} From the above discussion, for $G$ a cyclotomic graph, $R_G$ has all its
roots of modulus $1$ and so is a cyclotomic polynomial, by
Kronecker's Theorem.

We now take $G$ to be a Salem
graph, with index
$\lambda=\lambda(G)$. We can construct its reciprocal polynomial,
$R_G$, and $\tau=\tau(G)$ is a root of this. Moreover, $\lambda$
is the only root of $P_G$ that is greater than 2, so that apart
from $\lambda$ and possibly $-\lambda$ all the roots of $P_G$ lie
in the real interval $[-2,2\,]$. As noted above, such roots of
$P_G$ correspond to roots of $R_G$ that have modulus 1; $\la$
(respectively the pair $\pm\lambda$) corresponds to the pair of
real roots $\tau$, $1/\tau$, with $\tau>1$. The minimal polynomial
of $\tau$, call it $m_\tau$, is a factor of $R_G$. Its roots include
$\tau$ and $1/\tau$. Were  $\lambda$ (respectively $\la^2$) to be
a rational integer---cases excluded in the definition---then these
would be the only roots of $m_\tau$, and $\tau$ would be a
reciprocal quadratic Pisot number. As this is not the case,
$m_\tau$ has at least one root with modulus 1, and exactly one
root ($\tau$) with modulus greater than 1, so $\tau$ is a Salem
number.
\end{proof}

Many of the results that follow are most readily stated using
Salem graphs, although our real interest is only in nontrivial Salem
graphs. It is an easy matter, however, to check from the definition whether or not a
particular Salem graph is trivial.

While we are able in Section \ref{S-trees} to describe all Salem
trees, we are not at present able to do the same for Salem graphs.
However, the following result
  greatly restricts the kinds of
graphs that can be Salem graphs. It is an
essential ingredient in the proof of Theorem \ref{T-closed}.

\begin{proposition} \label{P-maxdeg3} Let $G$ be a connected graph having index  $\la>2$ and
second largest eigenvalue at most $2$. Then

\begin{enumerate}
\item[(a)] The vertices $V(G)$
 of $G$ can be partitioned as $V(G)=M\cup A\cup H$,
 in such a way that
 \begin{itemize}
 \item The induced subgraph $G|_{M}$ is one of the $18$ graphs
 of \cite[Theorem 2.3]{CR} minimal with respect to the property of
 having index greater than $2$; it has only one
 eigenvalue
 greater than $2$;
 \item The set $A$  consists of all vertices of $G-M$ adjacent in $G$ to
 some vertex of $M$;
 \item The induced subgraph $G|_H$ is cyclotomic.
 \end{itemize}

\item[(b)]  $G$ has  at most $B:=10(3\la^4 + \la^2 + 1)$ vertices
 of degree greater
than $2$, and at most $\la^2 B$ vertices of degree $1$.
\end{enumerate}
 \end{proposition}

 \begin{proof}
  Such a graph $G$  has a minimal
 vertex-deleted induced subgraph $G|_M$ with index  greater than $2$, given by \cite[Theorem 2.3]{CR};
 $G|_M$ can be one of $18$ graphs, each with at most $10$ vertices.
Note that $G|_M$ has only one eigenvalue greater than $2$, as when a
vertex is removed from $G|_M$ the resulting graph has, by
minimality, index at most $2$. Hence, by Lemma \ref{L-interlacing},
$G|_M$ cannot have more than one eigenvalue greater than $2$.

Now let $A$ be the set of vertices in $V(G)-M$ adjacent in $G$ to a
vertex
 of $M$.
 Then,  by interlacing, the induced subgraph  $G'$ on $V(G)-A$ has at most one eigenvalue greater than $2$, which must be the
index of $G|_M$. Hence the other components of $G'$ must together
form a cyclotomic graph, $H$
 say. By definition, there are no edges in $G$ having one endvertex in $M$ and the other
in $H$.

As the index of $G$ is $\la$, the maximum degree of a vertex of
$G$ is  bounded  by $\la^2$, by Lemma \ref{L-maxdeg}. Applying
this to the
 vertices of $M$, we see
that there are at most $10\la^2$ edges with one endvertex in
$M$ and the other in $A$. Thus the size $\# A$ of $A$ is at
most $10\la^2$.
 Now, applying the degree
bound $\la^2$ to the vertices of $A$, we similarly get the upper
bound $\la^2\#A $ for the number of edges with one endvertex in $A$
and the other in $H$. These edges are adjacent to at most
$\la^2\#A $ vertices in $H$ of degree $> 2$ in  $G$.
Now every
connected cyclotomic graph contains at most two vertices of
degree greater than $2$ (in fact only the type $\tilde D_n$, 
as in Figure \ref{F-except}, having two). Also, since every
connected component of $H$ has at least one such edge incident in
it, the number of such components is at most $\la^2\#A$. This gives
at most another $2\la^2\#A$ vertices of degree $> 2$ in $H$ that are
not adjacent to a vertex of $A$. Adding up, we see that the total
number of vertices of degree $> 2$ is at most $\#M + \#A+ \la^2\#A\
+ 2\la^2\#A\ \leq 10(3\la^4 + \la^2 + 1)$.

To bound the number of vertices of degree $1$, we associate to each such vertex
the nearest (in the obvious sense) vertex of degree greater than $2$, and then use the fact that these latter vertices have degree at most $\la(G)^2$, by Lemma \ref{L-maxdeg}.
\end{proof}

On the positive side, the next results enable us to construct many Salem graphs.
Our first result does this for bipartite Salem graphs.

\begin{theorem}\label{T-graphs}
\begin{enumerate}
\item[(a)] Suppose that $G$ is a noncyclotomic bipartite graph and such that
the induced subgraph
on $V(G)-\{v\}$ is cyclotomic. Then $G$ is a Salem graph.

\item[(b)] Suppose that $G$ is a noncyclotomic bipartite graph,  with the property
that for each minimal induced subgraph $M$ of $G$ the complementary
induced subgraph $G|_{V(G)-V(M)}$ is cyclotomic. Then $G$ is
a Salem graph.

 \end{enumerate}
\end{theorem}
Here the ``minimal'' graph $M$ is as in  Proposition
\ref{P-maxdeg3}: a minimal vertex-deleted subgraph with index
greater than $2$.

We can use  part (a) of the theorem to construct Salem graphs. Take
a forest of cyclotomic bipartite graphs (that is, any graph of Lemma
\ref{L-except} except an odd cycle $\tilde A_{2n}$), and colour the
vertices black or red, with adjacent vertices differently coloured.
Join some (as few or as many as you like) of the black vertices to a
new red vertex. Of course, one may as well take enough such edges to
make $G$ connected. This construction gives the most general
bipartite, connected graph such that removing the vertex $v$
produces a graph with all eigenvalues in $[-2,2]$. This result is an
extension of Theorem \ref{T-trees}(a) below, which is for trees.
Theorem \ref{T-trees}(b) gives a construction for more Salem trees.

In 2001 Piroska Lakatos \cite{L2} proved a special case of Theorem
\ref{T-graphs} where the components $G|_{V(G)-\{v\}}$ consisted of
paths, joined in $G$ at one or both endvertices to $v$.

\begin{proof}
The proof of (a) is immediate from Lemma \ref{L-interlacing}.

Part (b) comes straight from a result of D. Powers---see \cite[p.
456]{CS}. This states that if the vertices of a graph $G$ are
partitioned as $V(G)=V_1\cup V_2$ with $ G|_{V_i}(i=1,2)$ having
indices $\la^{(i)}(i=1,2)$, then the second-largest eigenvalue of
$G$ is at most $\max_{V_1\cup V_2=V(G)} \min(\la^{(1)},\la^{(2)})$.
It is clear that we may restrict consideration to $G|_{V_1}$ that
are minimal, which gives the result.
\end{proof}

The next result gives a construction for some nonbipartite Salem graphs.

\begin{theorem}\label{T-nonbip-graphs}
 Suppose that $G$ is a noncyclotomic nonbipartite graph containing a vertex $v$ such that
the induced subgraph
on $V(G)-\{v\}$ is cyclotomic. Suppose also that $G$ is a line graph. Then $G$ is a Salem graph.
\end{theorem}

\begin{proof}
Recall that a line graph $L$ is obtained from another graph $H$ by defining the vertices
of $L$ to be the edges of $H$, with two vertices of $L$ adjacent if and only if the
corresponding edges of $H$ are incident at a common vertex of $H$.
It is known that line graphs have least eigenvalue at least $-2$ (\cite[Chapter 12]{GR}).
The proof of this follows easily from the fact that the adjacency matrix $A$ of $L$ is given
by $A+2I=B^TB$, where $B$ is the incidence matrix of $H$.
Further, by Lemma \ref{L-interlacing}, $G$ has one eigenvalue $\la(G)>2$.
\end{proof}

 To use this result constructively, first note that all paths and cycles are line graphs, as well as being cyclotomic. Then take any graph $H$ consisting of one or two connected components, each of which is a path or cycle, and add to $H$ an extra edge joining any two distinct nonadjacent vertices. Then the line graph of this augmented graph,
if not again cyclotomic, will be a nonbipartite Salem
graph.

\section{Lemmas on reciprocal polynomials of graphs}\label{S-recpolys}
For the proof of Theorem \ref{T-closed}, we shall need to consider
special families of graphs, obtained by adding paths to a graph.
Here we establish the general structure of the reciprocal polynomials of
such families, and show how in certain cases one can retrieve
a Pisot number from a sequence of graph Salem numbers.

Throughout this section, reciprocal polynomials will be written as
functions of a variable $z$, and we conveniently treat the
bipartite and nonbipartite cases together by writing $y=\sqrt{z}$
if the graph is bipartite, and $y=z$ otherwise.

\begin{lemma}\label{L-1path}
Let $G$ be a graph with a distinguished vertex $v$. For each
$m\ge0$, let $G_m$ be the graph obtained by attaching one endvertex
of an $m$-vertex path to the vertex $v$ (so $G_m$ has $m$ more
vertices than $G$).

Let $R_m(z)$ be the reciprocal polynomial of $G_m$.
Then for $m\geq 2$ we have
\[
(y^2-1)R_m(z)=y^{2m}P(z)-P^*(z)\,,
\]
for some monic integer polynomial $P(z)$
that depends on $G$ and $v$ but not on $m$,
and with $P^{*}(z)=z^{\deg P}P(1/z)$.
\end{lemma}

\begin{proof}
Let $\pchi_m(\lambda)$ be the characteristic polynomial of $G_m$.
Then expanding this determinant along the row corresponding to the
vertex at the ``loose'' endvertex of the attached path (that which is not $v$)
we get (for $m\geq 2$)
\[
\pchi_m = \lambda\pchi_{m-1} - \pchi_{m-2}\,.
\]
Recognising this as a Chebyshev recurrence, or using induction, we
get (on replacing $\lambda$ by $y+1/y$ and multiplying through by
the appropriate power of $y$)
\[
R_m(z)=\frac{y^{2(t+1)}-1}{y^2-1}R_{m-t}(z)
-\frac{y^{2t}-1}{y^2-1}y^2R_{m-t-1}(z)
\]
for any $t$ between 1 and $m-1$.
Taking $t=m-1$ gives
\[
R_m(z)=\frac{y^{2m}-1}{y^2-1}R_1(z)-\frac{y^{2(m-1)}-1}{y^2-1}y^2R_0(z)\,.
\]
Putting $P(z)=R_1(z)-R_0(z)$ we are done.
\end{proof}

An easy induction extends this lemma to deal with any number of
added pendant paths.

\begin{lemma}\label{L-kpaths}
Let $G$ be a graph, and $(v_1,\ldots,v_k)$ a list of (not
necessarily distinct) vertices of $G$. Let $G_{m_1,\ldots,m_k}$ be
the graph obtained by attaching one endvertex of a new
$m_i$-vertex path to vertex $v_i$ (so $G_{m_1,\ldots,m_k}$ has
$m_1+\ldots+m_k$ more vertices than $G$).
Let $R_{m_1,\ldots,m_k}(z)$ be the reciprocal polynomial of
$G_{m_1,\ldots,m_k}$.
Then if all the $m_i$ are $\geq 2$ we have
\[
(y^2-1)^kR_{m_1,\ldots,m_k}(z)=
\sum_{\epsilon_1,\ldots,\epsilon_k\in\{0,1\}}
y^{2\sum\epsilon_im_i}
P_{(\epsilon_1,\ldots,\epsilon_k)}(z)\,,
\]
for some integer polynomials $P_{(\epsilon_1,\ldots,\epsilon_k)}(z)$
that depend on $G$ and $(v_1,\ldots,v_k)$ but not on $m_1,\ldots,m_k$.
\end{lemma}

With notation
as in the Lemma,
we refer to $P_{(1,\ldots,1)}$ as the {\em leading polynomial}
of $R_{m_1,\ldots,m_k}$.
Given any $\epsilon>0$, if we then take
all the $m_i$ large enough, the number of zeros of $R_{m_1,\ldots,m_k}$
in the region $|z|\geq 1+\epsilon$ is equal to the number of
zeros of its leading polynomial in that region.

\begin{lemma}\label{L-1pathpisot}
Suppose that $G$ is connected and that $G_m$ (as in Lemma
\ref{L-1path}) is a Salem graph for all sufficiently
large $m$. Then $G_m$ is a nontrivial Salem graph for all sufficiently
large $m$.
  Furthermore $P(z)$,
the leading polynomial of $R_m$, is a product of a Pisot polynomial
(with Pisot number $\theta$ as its root, say), a power of $z$, and
perhaps a cyclotomic polynomial. Moreover the Salem numbers $\tau_m:=\tau(G_m)$
converge to $\theta$ as $m\rightarrow\infty$.
\end{lemma}

\begin{proof}
Preserving the notation of Lemma \ref{L-1path}, we have
$(y^2-1)R_m(z) = y^{2m}P(z)-P^*(z)$. We suppose that $m$ is large
enough that the only roots of $R_m(z)$ are
$\tau_m$, its conjugates, and perhaps some roots of unity. By
Lemma \ref{L-addedge}, the $\tau_m$ are strictly increasing, so in
particular they have modulus $\geq 1+\epsilon$ for all sufficiently
small positive $\epsilon$. From the remark preceding this Lemma,
we deduce that $P(z)$ has exactly one root outside the closed unit
disc, $\theta$ say. Applying Rouch\'e's Theorem on the boundary of
an arbitrarily small disc centred on $\theta$ we deduce that, for
all large enough $m$, $R_m$ and $P$ have the same number of zeros
(namely one) within that disc, and hence $\tau_m\rightarrow\theta$
as $m\rightarrow\infty$. Since the eigenvalues of trivial Salem graphs form a discrete set,
we can discard the at most finite number of trivial Salem graphs in
our sequence, and so assume that all our Salem graphs are
nontrivial, so that the $\tau_m$ are Salem numbers.

It remains to prove that $\theta$ is a Pisot number.
The only alternative would be that $\theta$ is
a Salem number.
But then $\theta$ would also be a root of $P^*(z)$,
so would be a root of $R_m(z)$ for all $m$,
giving $\tau_m=\theta$ for all $m$.
This contradicts the fact that the $\tau_m$
are strictly increasing as $m$ increases.
\end{proof}

It is interesting to note that the Pisot number $\theta$ in Lemma
\ref{L-1pathpisot} cannot be a reciprocal quadratic Pisot number,
the proof showing that it is not conjugate to $1/\theta$.

\begin{corollary}\label{C-kpaths}
With notation as in Lemma \ref{L-kpaths}, suppose further that $G$
is connected and that $G_{m_1,\ldots,m_k}$ is a  Salem
graph for all sufficiently large $m_1,\ldots,m_k$.
Then $G_{m_1,\ldots,m_k}$ is a nontrivial Salem graph for but finitely many $(m_1,\ldots,m_k)$. Furthermore
$P_{(1,\ldots,1)}(z)$, the leading polynomial of
$R_{m_1,\ldots,m_k}(z)$, is the product of the minimal polynomial of
some Pisot number ($\theta$, say), a power of $z$, and perhaps a
cyclotomic polynomial.

Moreover, if we all let the $m_i$ tend to infinity in any
manner (one at a time, in bunches, or all together, perhaps
at varying rates),
the Salem numbers $\tau_{m_1,\ldots,m_k}=\tau(G_{m_1,\ldots,m_k})$
tend to $\theta$.
\end{corollary}

\begin{proof}
Throughout we suppose that the $m_i$ are all sufficiently large that
all the graphs under consideration are  Salem graphs. As in the proof of the previous lemma,
we may assume that these are nontrivial, so that the $\tau_{m_1,\ldots,m_k}$ are Salem numbers.
Fixing $m_2,\ldots, m_k$ (all large enough), and letting
$m_1\rightarrow\infty$, we apply Lemma \ref{L-1pathpisot} to deduce
that $\tau_{m_1,\ldots,m_k}$ tends to a Pisot number, say
$\theta_{\infty,m_2,\ldots,m_k}$, that is a root of
\[
\sum_{\epsilon_2,\ldots,\epsilon_k\in\{0,1\}}
y^{2\sum_{i\ge2}\epsilon_i m_i}
P_{(1,\epsilon_2,\ldots,\epsilon_k)}(z)\,.
\]

Now we let $m_2\rightarrow\infty$, and we get a sequence of
Pisot numbers that converge to the unique root of
\[
\sum_{\epsilon_3,\ldots,\epsilon_k\in\{0,1\}}
y^{2\sum_{i\ge3}\epsilon_i m_i}
P_{(1,1,\epsilon_3,\ldots,\epsilon_k)}(z)\,.
\]
outside the closed unit disc.
Since the set of Pisot numbers is closed (\cite{Sa}), this
number, $\theta_{\infty,\infty,m_3,\ldots,m_k}$,
must be a Pisot number.

Similarly we let the remaining $m_i\rightarrow\infty$,
producing a Pisot number $\theta=\theta_{\infty,\ldots,\infty}$
that is the unique root of $P_{(1,1,\ldots,1)}$
outside the closed unit disc.
Hence $P_{(1,1,\ldots,1)}$ has the desired form.

Finally we note that in whatever manner the $m_i$ tend to infinity,
$P_{(1,1,\ldots,1)}$ eventually dominates outside the unit circle,
and a Rouch\'e argument near $\theta$ shows that the Salem numbers
converge to $\theta$.
\end{proof}

\begin{lemma}\label{L-splitends}
Let $G$ be a graph with two (perhaps equal) distinguished vertices
$v_1$ and $v_2$. Let $G^{(m_1,m_2)}$ be the graph obtained by
identifying the endvertices of a  new $(m_1+m_2+3)$-vertex
path with vertices $v_1$ and $v_2$
(so that $G^{(m_1,m_2)}$ has $m_1+m_2+1$ more vertices than $G$).
Let $R^{(m_1,m_2)}$ be the reciprocal polynomial of $G^{(m_1,m_2)}$.

Removing the appropriate vertex ($w$ say) from the new path, we
get the graph $G_{m_1,m_2}$ (in the notation of Lemma
\ref{L-kpaths}), with reciprocal polynomial $R_{m_1,m_2}$.

Then
\[
R^{(m_1,m_2)}=(y^2-1)R_{m_1,m_2}(z)+Q_{m_1,m_2}(z)\,,
\]
where $Q_{m_1,m_2}$ has  negligible degree compared to
$R_{m_1,m_2}$, in the sense that
$$\deg(Q_{m_1,m_2})/\deg(R_{m_1,m_2})\rightarrow0$$ as
$\min(m_1,m_2)\rightarrow\infty$.
\end{lemma}

With the natural extension of our previous notion of
a leading polynomial, this Lemma implies that
$R^{(m_1,m_2)}$ has the same leading polynomial
as $R_{m_1,m_2}$.

\begin{proof}
Expanding $\pchi_{G^{(m_1,m_2)}}=\det(\lambda I - {\rm adjacency\
matrix\ of\ }G^{(m_1,m_2)})$ along the row corresponding to the
vertex $w$, we get
\[
\pchi_{G^{(m_1,m_2)}}=\lambda \pchi_{G_{m_1,m_2}}
-\pchi_{G_{m_1-1,m_2}}-\pchi_{G_{m_1,m_2-1}} +Q_1(\lambda)\,,
\]
where $Q_1$, and also $Q_2,Q_3,Q_4$ below, have negligible degree
compared to the other polynomials in the equation where they
appear.

Substituting $\lambda = y+1/y$
and multiplying by the appropriate power of $y$ gives
\[
R^{(m_1,m_2)} = (y^2+1)R_{m_1,m_2}(z)-y^2 R_{m_1-1,m_2}(z) -y^2
R_{m_1,m_2-1}(z) + Q_2(z)\,.
\]

Applying Lemma \ref{L-kpaths} for the case $k=2$, we get
\begin{align*}
(y^2-1)^2 R^{(m_1,m_2)}(z) =& P_{(1,1)}(z)
\left\{
(y^2+1)y^{2(m_1+m_2)}
-y^{2+2(m_1-1+m_2)}
-y^{2+2(m_1+m_2-1)}
\right\}
+Q_3(z)\\
=&y^{2(m_1+m_2)}(y^2-1)P_{(1,1)}(z)+Q_3(z)\,.
\end{align*}

Comparing with
\[
(y^2-1)^2R_{m_1,m_2}(z) = y^{2(m_1+m_2)}P_{(1,1)}(z)+Q_4(z)\,,
\]
 we get the advertised result.
\end{proof}

\section{proof of Theorem \ref{T-closed}}\label{S-proof1}

\begin{proof} Consider an infinite sequence of nontrivial Salem graphs $G$, for which the Salem numbers
$\tau(G)$ tend to a limit.
We are
interested in limit points of the set $T_{\rm graph}$, so we may
suppose, by moving to a subsequence, that our sequence has no
constant subsequence; moreover we can suppose that the graphs are
either all bipartite, or all nonbipartite. Indeed we shall suppose
that they are all nonbipartite, and leave the trivial modifications
for the bipartite case to the reader. These Salem numbers are
bounded above, and hence so are the indices of their graphs. Hence
Proposition \ref{P-maxdeg3} gives an upper bound on the number of
vertices of degree not equal to $2$ of these Salem graphs, and Lemma
\ref{L-maxdeg} gives an upper bound on the degrees of vertices that
each such graph can have.
 Now,  the set
of all multigraphs with at most $B_1$ vertices each of which is of
degree at most $B_2$ is finite. Thus, on associating to each Salem
graph in the sequence the multigraph with no vertices of degree $2$
having that Salem graph as a subdivision (that is, placing extra
vertices of degree $2$ along  edges of the multigraph retrieves the
Salem graph),
 we obtain only finitely many different
 multigraphs. Hence,
by replacing the sequence of Salem graphs by a subsequence, if necessary, we can assume that all Salem graphs in the
sequence are associated to the same multigraph, $M$ say. Now label the edges of $M$ by $e_1\cc e_m$ say. Each edge
$e_j$ corresponds to a path, of length $\ell_{j,n}$ say,
on the $n$-th Salem graph of the sequence, joining two vertices of degree not equal to $2$.

Now consider the sequence $\{\ell_{1,n}\}$. If it is bounded,
it has an infinite constant subsequence. Otherwise, it has a subsequence
tending monotonically to infinity. Hence, on taking a suitable
subsequence, we can assume that $\{\ell_{1,n}\}$ has one or other of these properties. Furthermore, since any infinite
subsequence of a sequence having one
of these properties inherits that property, we can take further infinite
subsequences without losing that property. Thus we do the same successively
for $\{\ell_{2,n}\}$, then $\{\ell_{3,n}\},\{\ell_{4,n}\},\dots,\{\ell_{m,n}\}$.
The effect is that we can assume that every sequence $\{\ell_{j,n}\}$ is
either constant or tends to infinity monotonically. Those that are constant
can simply be incorporated into $M$ (now allowing it to have vertices of degree $2$), so that we can in fact assume that
they all tend to infinity monotonically.

Let us suppose that our sequence of Salem graphs, $\{G_r\}$ has $s$
increasingly-subdivided internal edges, and $t$ pendant-increasing
edges. Form another set of graphs by removing a vertex from the
middle (or near middle) of each increasingly-subdivided edge of each
$G_r$, leaving $2s+t$ pendant-increasing edges.
We shall use $K_r$ to denote a
graph in this sequence,
 with $n_1$, \ldots, $n_{2s+t}$ for the lengths of its
pendant-increasing edges.

Claim: for {\em any} sufficiently large $n_1$, \ldots, $n_{2s+t}$,
we have a Salem graph. For (i) we soon exclude all cyclotomic graphs
from the list given in Section \ref{S-quotients}; and (ii) we can
never have more than one eigenvalue that is $>2$, otherwise,
 by adding vertices to reach one of our $G_r$ we would find a
Salem graph with more than one eigenvalue $>2$, using Lemma
\ref{L-interlacing}; and (iii) we can never have an eigenvalue
that is $<-2$, by similar reasoning.

Now we apply Corollary \ref{C-kpaths} to deduce that the limit of
our sequence of Salem numbers coming from the $K_r$ is a Pisot
number. (Note that $K_r$ need not be connected. All but one
component will be cyclotomic, and the noncyclotomic component
produces our Pisot number (the others merely contribute cyclotomic
factors to the leading polynomial).) Finally, by Lemma
\ref{L-splitends} this limiting Pisot number is also the limit of
the original sequence of Salem numbers.

The last sentence of Theorem \ref{T-closed} follows immediately.
\end{proof}

Examining the proof, we see that the number $m$ of lengthening
paths attached to the noncyclotomic growing component tells us
that the limiting Pisot number is in the $m$-th
 derived set of $T_{\gra}$, and so in the $(m-1)$-th
 derived set of $S_{\gra}$.

\section{Cyclotomic rooted trees}\label{S-quotients}

If $T$ is a rooted  tree, by which we of course mean a  tree
with a distinguished vertex $r$ say,  its {\em root},
 then $T'$ will
denote the rooted forest (set of rooted trees) $T-\{r\}$, the root
of each tree in $T'$ being its vertex  that is adjacent
(in $T$) to $r$.

 The {\it  quotient} of a rooted tree is the rational function
 $\pq_T=\frac{\prod_i R_i(z)}{R_{T}(z)}$, where $R_T$ is the
reciprocal polynomial of the tree, and the $R_i$ are the
reciprocal polynomials of its rooted subtrees, the trees
of $T'$.
 We define the {\em $\nu$-value} $\nu(T)$
of a tree $T$ to be $\pq_{T}(1)$, allowing $\nu(T)=\infty$ if
$\pq_T$ has a pole at 1. Note that by Lemma \ref{L-interlacing}, all
zeros and poles of  $(z-1)\pq_T$  are simple. The poles correspond
to a subset of the distinct eigenvalues of $T$ via $\lambda =
\sqrt{z}+1/\sqrt{z}$.

In this section we use Lemma \ref{L-except} to list all rooted
cyclotomic trees, along with their quotients and $\nu$-values. These
will be used in the following section (Theorem \ref{T-trees}) to
show how to construct all Salem trees.

 In our list, each entry for a tree $T$ contains
the following: a name for $T$, based on Coxeter graph notation;
a picture of $T$, with the root circled;  its quotient
$\pq_T(z)$ and $\nu$-values $\nu(T)=\pq_T(1)$.
 Here
$\Phi_n=\Phi_n(z)$  is the $n$-th cyclotomic polynomial.

First, the   rooted trees that are proper subtrees of $\tilde
E_6,\tilde E_7$ or $\tilde E_8$, but  not subtrees of any $\tilde D_n$:

\begin{picture}(350,45)
\put(0,25){$E_6(1)$} \put(35,20){\circle*{5}}
\put(35,20){\circle{10}} \put(35,20){\line(1,0){60}}
\put(50,20){\circle*{5}} \put(65,20){\circle*{5}}
\put(65,20){\line(0,1){15}} \put(65,35){\circle*{5}}
\put(80,20){\circle*{5}} \put(95,20){\circle*{5}}
\put(0,0){$\frac{\Phi_2\Phi_8}{\Phi_3\Phi_{12}}
=\frac{(z+1)(z^4+1)}{(z^2+z+1)(z^4-z^2+1)}, \nu=\frac{4}{3}$}
\put(220,25){$E_6(2)$} \put(255,20){\circle*{5}}
\put(255,20){\line(1,0){60}} \put(270,20){\circle*{5}}
\put(270,20){\circle{10}} \put(285,20){\circle*{5}}
\put(285,20){\line(0,1){15}} \put(285,35){\circle*{5}}
\put(300,20){\circle*{5}} \put(315,20){\circle*{5}}
\put(220,0){$\frac{\Phi_2\Phi_5}{\Phi_3\Phi_{12}}
=\frac{(z+1)(z^4+z^3+z^2+z+1)}{(z^2+z+1)(z^4-z^2+1)},
\nu=\frac{10}{3}$}
\end{picture}
\bigskip

\begin{picture}(350,40)
\put(0,25){$E_6(3)$} \put(35,20){\circle*{5}}
\put(65,20){\circle{10}} \put(35,20){\line(1,0){60}}
\put(50,20){\circle*{5}} \put(65,20){\circle*{5}}
\put(65,20){\line(0,1){15}} \put(65,35){\circle*{5}}
\put(80,20){\circle*{5}} \put(95,20){\circle*{5}}
\put(0,0){$\frac{\Phi_2\Phi_3}{\Phi_{12}}
=\frac{(z+1)(z^2+z+1)}{z^4-z^2+1}, \nu=6$}
 \put(220,25){$E_6(4)$}
\put(255,20){\circle*{5}} \put(255,20){\line(1,0){60}}
\put(270,20){\circle*{5}} \put(285,35){\circle{10}}
\put(285,20){\circle*{5}} \put(285,20){\line(0,1){15}}
\put(285,35){\circle*{5}} \put(300,20){\circle*{5}}
\put(315,20){\circle*{5}}
\put(220,0){$\frac{\Phi_2\Phi_6}{\Phi_{12}}
=\frac{(z+1)(z^2-z+1)}{z^4-z^2+1}, \nu=2$}
\end{picture}
\bigskip

\begin{picture}(350,45)
\put(0,25){$E_7(1)$} \put(35,20){\circle*{5}}
\put(35,20){\circle{10}} \put(35,20){\line(1,0){75}}
\put(50,20){\circle*{5}} \put(65,20){\circle*{5}}
\put(65,20){\line(0,1){15}} \put(65,35){\circle*{5}}
\put(80,20){\circle*{5}} \put(95,20){\circle*{5}}
\put(110,20){\circle*{5}}
\put(0,0){$\frac{\Phi_2\Phi_{10}}{\Phi_{18}}
=\frac{(z+1)(z^4-z^3+z^2-z+1)}{z^6-z^3+1}, \nu=2$}
\put(220,25){$E_7(2)$} \put(255,20){\circle*{5}}
\put(255,20){\line(1,0){75}} \put(270,20){\circle*{5}}
\put(270,20){\circle{10}} \put(285,20){\circle*{5}}
\put(285,20){\line(0,1){15}} \put(285,35){\circle*{5}}
\put(300,20){\circle*{5}} \put(315,20){\circle*{5}}
\put(330,20){\circle*{5}}
\put(220,0){$\frac{\Phi_2\Phi_3\Phi_6}{\Phi_{18}}
=\frac{(z+1)(z^2+z+1)(z^2-z+1)}{z^6-z^3+1}, \nu=6$}
\end{picture}
\bigskip

\begin{picture}(350,45)
\put(0,25){$E_7(3)$} \put(35,20){\circle*{5}}
\put(65,20){\circle{10}} \put(35,20){\line(1,0){75}}
\put(50,20){\circle*{5}} \put(65,20){\circle*{5}}
\put(65,20){\line(0,1){15}} \put(65,35){\circle*{5}}
\put(80,20){\circle*{5}} \put(95,20){\circle*{5}}
\put(110,20){\circle*{5}}
\put(0,0){$\frac{\Phi_2\Phi_3\Phi_4}{\Phi_{18}}
=\frac{(z+1)(z^2+z+1)(z^2+1)}{z^6-z^3+1}, \nu=12$}
\put(220,25){$E_7(4)$} \put(255,20){\circle*{5}}
\put(255,20){\line(1,0){75}} \put(270,20){\circle*{5}}
\put(300,20){\circle{10}} \put(285,20){\circle*{5}}
\put(285,20){\line(0,1){15}} \put(285,35){\circle*{5}}
\put(300,20){\circle*{5}} \put(315,20){\circle*{5}}
\put(330,20){\circle*{5}}
\put(220,0){$\frac{\Phi_3\Phi_5}{\Phi_2\Phi_{18}}
=\frac{(z^2+z+1)(z^4+z^3+z^2+z+1)}{(z+1)(z^6-z^3+1)},
\nu=\frac{15}{2}$}
\end{picture}
\bigskip

\begin{picture}(350,45)
\put(0,25){$E_7(5)$} \put(35,20){\circle*{5}}
\put(95,20){\circle{10}} \put(35,20){\line(1,0){75}}
\put(50,20){\circle*{5}} \put(65,20){\circle*{5}}
\put(65,20){\line(0,1){15}} \put(65,35){\circle*{5}}
\put(80,20){\circle*{5}} \put(95,20){\circle*{5}}
\put(110,20){\circle*{5}}
\put(0,0){$\frac{\Phi_2\Phi_8}{\Phi_{18}}
=\frac{(z+1)(z^4+1)}{z^6-z^3+1}, \nu=4$}
 \put(220,25){$E_7(6)$}
\put(255,20){\circle*{5}} \put(255,20){\line(1,0){75}}
\put(270,20){\circle*{5}} \put(330,20){\circle{10}}
\put(285,20){\circle*{5}} \put(285,20){\line(0,1){15}}
\put(285,35){\circle*{5}} \put(300,20){\circle*{5}}
\put(315,20){\circle*{5}} \put(330,20){\circle*{5}}
\put(220,0){$\frac{\Phi_3\Phi_{12}}{\Phi_2\Phi_{18}}
=\frac{(z^2+z+1)(z^4-z^2+1)}{(z+1)(z^6-z^3+1)}, \nu=\frac{3}{2}$}
\end{picture}
\bigskip

\begin{picture}(350,45)
\put(0,25){$E_7(7)$} \put(35,20){\circle*{5}}
\put(65,35){\circle{10}} \put(35,20){\line(1,0){75}}
\put(50,20){\circle*{5}} \put(65,20){\circle*{5}}
\put(65,20){\line(0,1){15}} \put(65,35){\circle*{5}}
\put(80,20){\circle*{5}} \put(95,20){\circle*{5}}
\put(110,20){\circle*{5}} \put(0,0){$\frac{\Phi_7}{\Phi_2\Phi_{18}}
=\frac{z^6+z^5+z^4+z^3+z^2+z+1}{(z+1)(z^6-z^3+1)}, \nu=\frac{7}{2}$}
\end{picture}
\bigskip
%"" $\tilde E_6$ moved""

\begin{picture}(350,45)
\put(0,25){$E_8(1)$} \put(35,20){\circle*{5}}
\put(35,20){\circle{10}} \put(35,20){\line(1,0){90}}
\put(50,20){\circle*{5}} \put(65,20){\circle*{5}}
\put(65,20){\line(0,1){15}} \put(65,35){\circle*{5}}
\put(80,20){\circle*{5}} \put(95,20){\circle*{5}}
\put(110,20){\circle*{5}} \put(125,20){\circle*{5}}
\put(0,0){$\frac{\Phi_2\Phi_4\Phi_{12}}{\Phi_{30}}
=\frac{(z+1)(z^2+1)(z^4-z^2+1)}{z^8+z^7-z^5-z^4-z^3+z+1}, \nu=4$}
\put(220,25){$E_8(2)$} \put(255,20){\circle*{5}}
\put(255,20){\line(1,0){90}} \put(270,20){\circle*{5}}
\put(270,20){\circle{10}} \put(285,20){\circle*{5}}
\put(285,20){\line(0,1){15}} \put(285,35){\circle*{5}}
\put(300,20){\circle*{5}} \put(315,20){\circle*{5}}
\put(330,20){\circle*{5}} \put(345,20){\circle*{5}}
\put(220,0){$\frac{\Phi_2\Phi_7}{\Phi_{30}}
=\frac{(z+1)(z^6+z^5+z^4+z^3+z^2+z+1)}{z^8+z^7-z^5-z^4-z^3+z+1}, \nu=14$}

\end{picture}
\bigskip

\begin{picture}(350,45)
\put(0,25){$E_8(3)$} \put(35,20){\circle*{5}}
\put(65,20){\circle{10}} \put(35,20){\line(1,0){90}}
\put(50,20){\circle*{5}} \put(65,20){\circle*{5}}
\put(65,20){\line(0,1){15}} \put(65,35){\circle*{5}}
\put(80,20){\circle*{5}} \put(95,20){\circle*{5}}
\put(110,20){\circle*{5}} \put(125,20){\circle*{5}}
\put(-10,0){$\frac{\Phi_2\Phi_3\Phi_5}{\Phi_{30}}
=\frac{(z+1)(z^2+z+1)(z^4+z^3+z^2+z+1)}{z^8+z^7-z^5-z^4-z^3+z+1}, \nu=30$},

 \put(220,25){$E_8(4)$} \put(255,20){\circle*{5}}
\put(255,20){\line(1,0){90}} \put(270,20){\circle*{5}}
\put(300,20){\circle{10}} \put(285,20){\circle*{5}}
\put(285,20){\line(0,1){15}} \put(285,35){\circle*{5}}
\put(300,20){\circle*{5}} \put(315,20){\circle*{5}}
\put(330,20){\circle*{5}} \put(345,20){\circle*{5}}
\put(220,0){$\frac{\Phi_2\Phi_4\Phi_5}{\Phi_{30}}
=\frac{(z+1)(z^2+1)(z^4+z^3+z^2+z+1)}{z^8+z^7-z^5-z^4-z^3+z+1}, \nu=20$}

\end{picture}
\bigskip

\begin{picture}(350,45)
\put(0,25){$E_8(5)$} \put(35,20){\circle*{5}}
\put(95,20){\circle{10}} \put(35,20){\line(1,0){90}}
\put(50,20){\circle*{5}} \put(65,20){\circle*{5}}
\put(65,20){\line(0,1){15}} \put(65,35){\circle*{5}}
\put(80,20){\circle*{5}} \put(95,20){\circle*{5}}
\put(110,20){\circle*{5}} \put(125,20){\circle*{5}}
\put(0,0){$\frac{\Phi_2\Phi_3\Phi_8}{\Phi_{30}}
=\frac{(z+1)(z^2+z+1)(z^4+1)}{z^8+z^7-z^5-z^4-z^3+z+1}, \nu=12$}
\put(220,25){$E_8(6)$} \put(255,20){\circle*{5}}
\put(255,20){\line(1,0){90}} \put(270,20){\circle*{5}}
\put(330,20){\circle{10}} \put(285,20){\circle*{5}}
\put(285,20){\line(0,1){15}} \put(285,35){\circle*{5}}
\put(300,20){\circle*{5}} \put(315,20){\circle*{5}}
\put(330,20){\circle*{5}} \put(345,20){\circle*{5}}
\put(220,0){$\frac{\Phi_2\Phi_3\Phi_{12}}{\Phi_{30}}
=\frac{(z+1)(z^2+z+1)(z^4-z^2+1)}{z^8+z^7-z^5-z^4-z^3+z+1}, \nu=6$}
\end{picture}
\bigskip

\begin{picture}(350,45)
\put(0,25){$E_8(7)$} \put(35,20){\circle*{5}}
\put(125,20){\circle{10}} \put(35,20){\line(1,0){90}}
\put(50,20){\circle*{5}} \put(65,20){\circle*{5}}
\put(65,20){\line(0,1){15}} \put(65,35){\circle*{5}}
\put(80,20){\circle*{5}} \put(95,20){\circle*{5}}
\put(110,20){\circle*{5}} \put(125,20){\circle*{5}}
\put(0,0){$\frac{\Phi_2\Phi_{18}}{\Phi_{30}}
=\frac{(z+1)(z^6-z^3+1)}{z^8+z^7-z^5-z^4-z^3+z+1}, \nu=2$}
\put(220,25){$E_8(8)$} \put(255,20){\circle*{5}}
\put(255,20){\line(1,0){90}} \put(270,20){\circle*{5}}
\put(285,35){\circle{10}} \put(285,20){\circle*{5}}
\put(285,20){\line(0,1){15}} \put(285,35){\circle*{5}}
\put(300,20){\circle*{5}} \put(315,20){\circle*{5}}
\put(330,20){\circle*{5}} \put(345,20){\circle*{5}}
\put(220,0){$\frac{\Phi_2\Phi_4\Phi_8}{\Phi_{30}}
=\frac{(z+1)(z^2+1)(z^4+1)}{z^8+z^7-z^5-z^4-z^3+z+1}, \nu=8$}
\end{picture}
\bigskip

Next, the rooted versions of $\tilde E_6,\tilde E_7$ and $\tilde
E_8$. Note that all their quotients have a pole at $z=1$, so that
$\nu=\infty$ for all these trees.

\begin{picture}(350,60)
\put(0,25){$\tilde{E}_6(1)$} \put(35,20){\circle*{5}}
\put(35,20){\circle{10}} \put(35,20){\line(1,0){60}}
\put(50,20){\circle*{5}} \put(65,20){\circle*{5}}
\put(65,20){\line(0,1){30}} \put(65,35){\circle*{5}}
\put(65,50){\circle*{5}} \put(80,20){\circle*{5}}
\put(95,20){\circle*{5}}
\put(0,0){$\frac{\Phi_{12}}{\Phi_1^2\Phi_2\Phi_3}
=\frac{z^4-z^2+1}{(z-1)^2(z+1)(z^2+z+1)}$}
\put(220,25){$\tilde{E}_6(2)$} \put(255,20){\circle*{5}}
\put(255,20){\line(1,0){60}} \put(270,20){\circle*{5}}
\put(270,20){\circle{10}} \put(285,20){\circle*{5}}
\put(285,20){\line(0,1){30}} \put(285,35){\circle*{5}}
\put(285,50){\circle*{5}} \put(300,20){\circle*{5}}
\put(315,20){\circle*{5}}
\put(220,0){$\frac{\Phi_2\Phi_6}{\Phi_1^2\Phi_3}
=\frac{(z+1)(z^2-z+1)}{(z-1)^2(z^2+z+1)}$}
\end{picture}
\bigskip

\begin{picture}(350,60)
\put(0,25){$\tilde{E}_6(3)$} \put(35,20){\circle*{5}}
\put(65,20){\circle{10}} \put(35,20){\line(1,0){60}}
\put(50,20){\circle*{5}} \put(65,20){\circle*{5}}
\put(65,20){\line(0,1){30}} \put(65,35){\circle*{5}}
\put(65,50){\circle*{5}} \put(80,20){\circle*{5}}
\put(95,20){\circle*{5}} \put(0,0){$\frac{\Phi_3}{\Phi_1^2\Phi_2}
=\frac{z^2+z+1}{(z-1)^2(z+1)}$}
\end{picture}
\bigskip

\begin{picture}(350,45)
\put(0,25){$\tilde{E}_7(1)$} \put(35,20){\circle*{5}}
\put(35,20){\circle{10}} \put(35,20){\line(1,0){90}}
\put(50,20){\circle*{5}} \put(65,20){\circle*{5}}
\put(80,20){\line(0,1){15}} \put(80,35){\circle*{5}}
\put(80,20){\circle*{5}} \put(95,20){\circle*{5}}
\put(110,20){\circle*{5}} \put(125,20){\circle*{5}}
\put(0,0){$\frac{\Phi_{18}}{\Phi_1^2\Phi_2\Phi_3\Phi_4}
=\frac{z^6-z^3+1}{(z-1)^2(z+1)(z^2+z+1)(z^2+1)}$}
\put(220,25){$\tilde{E}_7(2)$} \put(255,20){\circle*{5}}
\put(255,20){\line(1,0){90}} \put(270,20){\circle*{5}}
\put(270,20){\circle{10}} \put(285,20){\circle*{5}}
\put(300,20){\line(0,1){15}} \put(300,35){\circle*{5}}
\put(300,20){\circle*{5}} \put(315,20){\circle*{5}}
\put(330,20){\circle*{5}} \put(345,20){\circle*{5}}
\put(220,0){$\frac{\Phi_2\Phi_{10}}{\Phi_1^2\Phi_3\Phi_4}
=\frac{(z+1)(z^4-z^3+z^2-z+1)}{(z-1)^2(z^2+z+1)(z^2+1)}$}
\end{picture}
\bigskip

\begin{picture}(350,45)
\put(0,25){$\tilde{E}_7(3)$} \put(35,20){\circle*{5}}
\put(65,20){\circle{10}} \put(35,20){\line(1,0){90}}
\put(50,20){\circle*{5}} \put(65,20){\circle*{5}}
\put(80,20){\line(0,1){15}} \put(80,35){\circle*{5}}
\put(80,20){\circle*{5}} \put(95,20){\circle*{5}}
\put(110,20){\circle*{5}} \put(125,20){\circle*{5}}
\put(-10,0){$\frac{\Phi_3\Phi_6}{\Phi_1^2\Phi_2\Phi_4}
=\frac{(z^2+z+1)(z^2-z+1)}{(z-1)^2(z+1)(z^2+1)}$}
\put(220,25){$\tilde{E}_7(4)$} \put(255,20){\circle*{5}}
\put(255,20){\line(1,0){90}} \put(270,20){\circle*{5}}
\put(300,20){\circle{10}} \put(285,20){\circle*{5}}
\put(300,20){\line(0,1){15}} \put(300,35){\circle*{5}}
\put(300,20){\circle*{5}} \put(315,20){\circle*{5}}
\put(330,20){\circle*{5}} \put(345,20){\circle*{5}}
\put(220,0){$\frac{\Phi_2\Phi_4}{\Phi_1^2\Phi_3}
=\frac{(z+1)(z^2+1)}{(z-1)^2(z^2+z+1)}$}
\end{picture}
\bigskip

\begin{picture}(350,45)
\put(0,25){$\tilde{E}_7(5)$} \put(35,20){\circle*{5}}
\put(80,35){\circle{10}} \put(35,20){\line(1,0){90}}
\put(50,20){\circle*{5}} \put(65,20){\circle*{5}}
\put(80,20){\line(0,1){15}} \put(80,35){\circle*{5}}
\put(80,20){\circle*{5}} \put(95,20){\circle*{5}}
\put(110,20){\circle*{5}} \put(125,20){\circle*{5}}
\put(0,0){$\frac{\Phi_8}{\Phi_1^2\Phi_2\Phi_3}
=\frac{z^4+1}{(z-1)^2(z+1)(z^2+z+1)}$}
\end{picture}
\bigskip

\begin{picture}(360,45)
\put(0,25){$\tilde E_8(1)$} \put(35,20){\circle*{5}}
\put(35,20){\circle{10}} \put(35,20){\line(1,0){105}}
\put(50,20){\circle*{5}} \put(65,20){\circle*{5}}
\put(65,20){\line(0,1){15}} \put(65,35){\circle*{5}}
\put(80,20){\circle*{5}} \put(95,20){\circle*{5}}
\put(110,20){\circle*{5}} \put(125,20){\circle*{5}}
\put(140,20){\circle*{5}}
\put(-10,0){$\frac{\Phi_2\Phi_{14}}{\Phi_1^2\Phi_3\Phi_5}
=\frac{(z+1)(z^6-z^5+z^4-z^3+z^2-z+1)}
{(z-1)^2(z^2+z+1)(z^4+z^3+z^2+z+1)}$}
 \put(220,25){$\tilde E_8(2)$}
\put(255,20){\circle*{5}} \put(255,20){\line(1,0){105}}
\put(270,20){\circle*{5}} \put(270,20){\circle{10}}
\put(285,20){\circle*{5}} \put(285,20){\line(0,1){15}}
\put(285,35){\circle*{5}} \put(300,20){\circle*{5}}
\put(315,20){\circle*{5}} \put(330,20){\circle*{5}}
\put(345,20){\circle*{5}} \put(360,20){\circle*{5}}
\put(220,0){$\frac{\Phi_2\Phi_4\Phi_8}{\Phi_1^2\Phi_3\Phi_5}
=\frac{(z+1)(z^2+1)(z^4+1)}{(z-1)^2(z^2+z+1)(z^4+z^3+z^2+z+1)}$}
\end{picture}
\bigskip

\begin{picture}(360,45)
\put(0,25){$\tilde E_8(3)$} \put(35,20){\circle*{5}}
\put(65,20){\circle{10}} \put(35,20){\line(1,0){105}}
\put(50,20){\circle*{5}} \put(65,20){\circle*{5}}
\put(65,20){\line(0,1){15}} \put(65,35){\circle*{5}}
\put(80,20){\circle*{5}} \put(95,20){\circle*{5}}
\put(110,20){\circle*{5}} \put(125,20){\circle*{5}}
\put(140,20){\circle*{5}}
\put(0,0){$\frac{\Phi_2\Phi_3\Phi_6}{\Phi_1^2\Phi_5}
=\frac{(z+1)(z^2+z+1)(z^2-z+1)}{(z-1)^2(z^4+z^3+z^2+z+1)}$}
 \put(220,25){$\tilde E_8(4)$} \put(255,20){\circle*{5}}
\put(255,20){\line(1,0){105}} \put(270,20){\circle*{5}}
\put(300,20){\circle{10}} \put(285,20){\circle*{5}}
\put(285,20){\line(0,1){15}} \put(285,35){\circle*{5}}
\put(300,20){\circle*{5}} \put(315,20){\circle*{5}}
\put(330,20){\circle*{5}} \put(345,20){\circle*{5}}
\put(360,20){\circle*{5}}
\put(220,0){$\frac{\Phi_5}{\Phi_1^2\Phi_2\Phi_3}
=\frac{z^4+z^3+z^2+z+1}{(z-1)^2(z+1)(z^2+z+1)}$}
\end{picture}
\bigskip

\begin{picture}(360,45)
\put(0,25){$\tilde E_8(5)$} \put(35,20){\circle*{5}}
\put(95,20){\circle{10}} \put(35,20){\line(1,0){105}}
\put(50,20){\circle*{5}} \put(65,20){\circle*{5}}
\put(65,20){\line(0,1){15}} \put(65,35){\circle*{5}}
\put(80,20){\circle*{5}} \put(95,20){\circle*{5}}
\put(110,20){\circle*{5}} \put(125,20){\circle*{5}}
\put(140,20){\circle*{5}}
\put(-10,0){$\frac{\Phi_2\Phi_4\Phi_8}{\Phi_1^2\Phi_3\Phi_5}
=\frac{(z+1)(z^2+1)(z^4+1)}{(z-1)^2(z^2+z+1)(z^4+z^3+z^2+z+1)}$}
 \put(220,25){$\tilde E_8(6)$} \put(255,20){\circle*{5}}
\put(255,20){\line(1,0){105}} \put(270,20){\circle*{5}}
\put(330,20){\circle{10}} \put(285,20){\circle*{5}}
\put(285,20){\line(0,1){15}} \put(285,35){\circle*{5}}
\put(300,20){\circle*{5}} \put(315,20){\circle*{5}}
\put(330,20){\circle*{5}} \put(345,20){\circle*{5}}
\put(360,20){\circle*{5}}
\put(220,0){$\frac{\Phi_3\Phi_{12}}{\Phi_1^2\Phi_2\Phi_5}
=\frac{(z^2+z+1)(z^4-z^2+1)}{(z-1)^2(z+1)(z^4+z^3+z^2+z+1)}$}
\end{picture}
\bigskip

\begin{picture}(360,45)
\put(0,25){$\tilde E_8(7)$} \put(35,20){\circle*{5}}
\put(125,20){\circle{10}} \put(35,20){\line(1,0){105}}
\put(50,20){\circle*{5}} \put(65,20){\circle*{5}}
\put(65,20){\line(0,1){15}} \put(65,35){\circle*{5}}
\put(80,20){\circle*{5}} \put(95,20){\circle*{5}}
\put(110,20){\circle*{5}} \put(125,20){\circle*{5}}
\put(140,20){\circle*{5}}
\put(-10,0){$\frac{\Phi_2\Phi_{18}}{\Phi_1^2\Phi_3\Phi_5}
=\frac{(z+1)(z^6-z^3+1)}{(z-1)^2(z^2+z+1)(z^4+z^3+z^2+z+1)}$}
 \put(220,25){$\tilde E_8(8)$} \put(255,20){\circle*{5}}
\put(255,20){\line(1,0){105}} \put(270,20){\circle*{5}}
\put(360,20){\circle{10}} \put(285,20){\circle*{5}}
\put(285,20){\line(0,1){15}} \put(285,35){\circle*{5}}
\put(300,20){\circle*{5}} \put(315,20){\circle*{5}}
\put(330,20){\circle*{5}} \put(345,20){\circle*{5}}
\put(360,20){\circle*{5}}
\put(220,0){$\frac{\Phi_{30}}{\Phi_1^2\Phi_2\Phi_3\Phi_5}
=\frac{z^8+z^7-z^5-z^4-z^3+z+1}
{(z-1)^2(z+1)(z^2+z+1)(z^4+z^3+z^2+z+1)}$}
\end{picture}
\bigskip

\begin{picture}(360,45)
\put(0,25){$\tilde E_8(9)$} \put(35,20){\circle*{5}}
\put(65,35){\circle{10}} \put(35,20){\line(1,0){105}}
\put(50,20){\circle*{5}} \put(65,20){\circle*{5}}
\put(65,20){\line(0,1){15}} \put(65,35){\circle*{5}}
\put(80,20){\circle*{5}} \put(95,20){\circle*{5}}
\put(110,20){\circle*{5}} \put(125,20){\circle*{5}}
\put(140,20){\circle*{5}}
\put(0,0){$\frac{\Phi_9}{\Phi_1^2\Phi_2\Phi_5}
=\frac{z^6+z^3+1}{(z-1)^2(z+1)(z^4+z^3+z^2+z+1)}$}
\end{picture}
\bigskip

Then, the five infinite families:

\begin{picture}(250,70)
\put(0,30){$A_n(a,b)$} \put(50,30){\circle*{5}}
\put(50,30){\line(1,0){10}} \put(65,30){\circle*{2}}
\put(70,30){\circle*{2}} \put(75,30){\circle*{2}}
\put(80,30){\line(1,0){10}} \put(90,30){\circle*{5}}
\put(90,30){\circle{10}} \put(90,30){\line(1,0){10}}
\put(105,30){\circle*{2}} \put(110,30){\circle*{2}}
\put(115,30){\circle*{2}} \put(120,30){\line(1,0){10}}
\put(130,30){\circle*{5}} \put(150,30){$n=a+b-1$ vertices, $1\leq
a\leq b$} \put(90,47){\vector(0,-1){10}} \put(30,50){($a$th from
left, $b$th from right)} \put(0,0){
$\frac{(z^a-1)(z^b-1)}{(z-1)(z^{a+b}-1)}, \nu=\frac{ab}{a+b}$}
\end{picture}
\bigskip

\begin{picture}(250,70)
\put(0,30){$D_n(a,b)$} \put(50,30){\circle*{5}}
\put(50,30){\line(1,0){10}} \put(65,30){\circle*{2}}
\put(70,30){\circle*{2}} \put(75,30){\circle*{2}}
\put(80,30){\line(1,0){10}} \put(90,30){\circle*{5}}
\put(90,30){\circle{10}} \put(90,30){\line(1,0){10}}
\put(105,30){\circle*{2}} \put(110,30){\circle*{2}}
\put(115,30){\circle*{2}} \put(120,30){\line(1,0){10}}
\put(130,30){\circle*{5}} \put(150,45){\circle*{5}}
\put(150,15){\circle*{5}} \put(130,30){\line(4,3){20}}
\put(130,30){\line(4,-3){20}} \put(150,30){$n=a+b$ vertices,
$a\geq 1$, $b\geq 2$} \put(90,47){\vector(0,-1){10}}
\put(30,50){($a$th from left, $b$th from right)} \put(0,0){
$\frac{(z^a-1)(z^{b-1}+1)}{(z-1)(z^{a+b-1}+1)}, \nu=a$}
\end{picture}
\bigskip

\begin{picture}(250,70)
\put(0,30){$D_n(0)$} \put(50,30){\circle*{5}}
\put(50,30){\line(1,0){10}} \put(70,30){\circle*{2}}
\put(80,30){\circle*{2}} \put(90,30){\circle*{2}}
\put(100,30){\circle*{2}} \put(110,30){\circle*{2}}
\put(120,30){\line(1,0){10}} \put(130,30){\circle*{5}}
\put(150,45){\circle*{5}} \put(150,45){\circle{10}}
\put(150,15){\circle*{5}} \put(130,30){\line(4,3){20}}
\put(130,30){\line(4,-3){20}} \put(150,30){$n$ vertices, $n\geq 5$
(note that $D_4(0)=D_4(1,3)$)} \put(0,0){
$\frac{z^n-1}{(z-1)(z+1)(z^{n-1}+1)}, \nu=\frac{n}{4}$}
\end{picture}
\bigskip

\begin{picture}(250,70)
\put(0,30){$\tilde{D}_n(a,b)$} \put(50,45){\circle*{5}}
\put(50,15){\circle*{5}} \put(70,30){\line(-4,3){20}}
\put(70,30){\line(-4,-3){20}} \put(70,30){\circle*{5}}
\put(70,30){\line(1,0){10}} \put(85,30){\circle*{2}}
\put(90,30){\circle*{2}} \put(95,30){\circle*{2}}
\put(100,30){\line(1,0){10}} \put(110,30){\circle*{5}}
\put(110,30){\circle{10}} \put(110,30){\line(1,0){10}}
\put(125,30){\circle*{2}} \put(130,30){\circle*{2}}
\put(135,30){\circle*{2}} \put(140,30){\line(1,0){10}}
\put(150,30){\circle*{5}} \put(170,45){\circle*{5}}
\put(170,15){\circle*{5}} \put(150,30){\line(4,3){20}}
\put(150,30){\line(4,-3){20}} \put(170,30){$n+1=a+b+1$ vertices,
$2\leq a\leq b$} \put(110,47){\vector(0,-1){10}} \put(50,50){($a$th
from left, $b$th from right)} \put(0,0){
$\frac{(z^{a-1}+1)(z^{b-1}+1)}{(z-1)(z^{a+b-2}-1)}, \nu=\infty$}
\end{picture}
\bigskip

\begin{picture}(250,70)
\put(0,30){$\tilde{D}_n(0)$} \put(50,45){\circle*{5}}
\put(50,15){\circle*{5}} \put(70,30){\line(-4,3){20}}
\put(70,30){\line(-4,-3){20}} \put(70,30){\circle*{5}}
\put(70,30){\line(1,0){10}} \put(90,30){\circle*{2}}
\put(100,30){\circle*{2}} \put(110,30){\circle*{2}}
\put(120,30){\circle*{2}} \put(130,30){\circle*{2}}
\put(140,30){\line(1,0){10}} \put(150,30){\circle*{5}}
\put(170,45){\circle*{5}} \put(170,45){\circle{10}}
\put(170,15){\circle*{5}} \put(150,30){\line(4,3){20}}
\put(150,30){\line(4,-3){20}} \put(170,30){$n+1$ vertices, $n\geq
4$} \put(0,0){ $\frac{z^{n-1}+1}{(z-1)(z+1)(z^{n-2}-1)},
\nu=\infty$}
\end{picture}
\bigskip

We also note in passing the rooted even cycles:

\begin{picture}(250,70)
\put(0,30){$\tilde{A}_{2n-1}$} \put(50,30){\circle*{5}}
\put(50,30){\circle{10}} \put(50,30){\line(1,0){10}}
\put(70,30){\circle*{2}} \put(80,30){\circle*{2}}
\put(90,30){\circle*{2}} \put(100,30){\circle*{2}}
\put(110,30){\circle*{2}} \put(120,30){\line(1,0){10}}
\put(130,30){\circle*{5}} \put(50,30){\line(0,1){20}}
\put(50,50){\circle*{5}} \put(50,50){\line(1,0){10}}
\put(70,50){\circle*{2}} \put(80,50){\circle*{2}}
\put(90,50){\circle*{2}} \put(100,50){\circle*{2}}
\put(110,50){\circle*{2}} \put(120,50){\line(1,0){10}}
\put(130,50){\circle*{5}} \put(130,30){\line(0,1){20}}
\put(150,30){$2n$ vertices, $n\geq 2$} \put(0,0){
$\frac{z^n+1}{(z-1)(z^n-1)}, \nu=\infty$}
\end{picture}
\bigskip

As we have seen in Theorem \ref{T-graphs}, these can be used for constructing bipartite Salem graphs,
but obviously not  Salem trees.
\bigskip

Note that, since $\tilde E_8(2)$ and $\tilde E_8(5)$ have the same quotient,
we can readily construct different Salem trees having the same quotient,
and hence corresponding to the same Salem number.

For each of the Salem quotients $S(z)$ catalogued above, we
observe in passing that $(z-1)S(z)$ is an interlacing quotient, as
defined in \cite{MS}. This is an easy consequence of the
Interlacing Theorem  (Lemma \ref{L-interlacing}).

\section{A complete description of Salem trees}\label{S-trees}
In this section we consider the case of those (of course
bipartite) Salem graphs defined by trees. As before, if $T$ is a
rooted tree, then $T'$ will denote the rooted forest obtained by
deleting the root $r$ of $T$, with the root of each subtree being
the vertex that (in $T$) is adjacent to $r$. The quotient of a
rooted forest is defined to be the sum of the quotients of its
rooted trees. For rooted trees $T_1$ and $T_2$, we define the
rooted tree $T_1+T_2$ to be the
 tree obtained by joining the roots of $T_1$ and $T_2$ by an edge, and making
  the root of $T_1$ its root.

\begin{lemma}\label{L-quotient}
\begin{enumerate}
\item[(i)] \cite [Corollary 4]{MRS} For a rooted tree $T$ with
rooted subtrees $T'=\{T_i\}$, its quotient $\pq_T$ is given
recursively by
$$
\pq_T=\frac{1}{z+1-z\pq_{T'}}=\frac{1}{z+1-z\sum_i\pq_{T_i}},
$$
with $q_\bullet=1/(z+1)$ for the single-vertex tree $\bullet$.
\item[(ii)] For the rooted tree $T_1+T_2$  we have
$$
\pq_{T_1+T_2}=\frac{\pq_{T_1}}{1-z\pq_{T_1}\pq_{T_2}}=
\frac{z+1-z\pq_{T'_2}}{(z+1-z\pq_{T'_1})(z+1-z\pq_{T'_2})-z}.
$$
\end{enumerate}
\end{lemma}

\begin{proof}[Proof of (ii)] Applying (i) to $T_1+T_2$ and then
to $T_1$ gives
\[
\pq_{T_1+T_2}=\frac{1}{z+1-z\pq_{T'_1}-z\pq_{T_2}}
=\frac{1}{1/\pq_{T_1}-z\pq_{T_2}}=\frac{\pq_{T_1}}{1-z\pq_{T_1}\pq_{T_2}}\,.
\]
Now applying (i) again to both $T_1$ and $T_2$ gives the alternative formula.
\end{proof}

Note that (i) implies that $\nu(T)=1/(2-\nu(T'))$, with
$\nu(T')=\sum_i\nu(T_i)$.
\bigskip

The next Theorem describes all Salem trees. For an alternative
approach to a generalisation of this topic, see Neumaier
\cite[Theorem 2.6]{Neu}.

\begin{theorem}\label{T-trees}
\begin{enumerate}
\item[(a)]
Suppose that $T$ is a rooted tree with $\nu(T')>2$,  for which the
forest $T'$ is a collection of cyclotomic trees. Then $T$ is a
Salem tree. (If $\nu(T')\leq 2$ then $T$ is again an cyclotomic
tree.)

\item[(b)] Suppose that $T_1$ and $T_2$ are Salem trees of type (a) with
$(\nu(T'_1)-2)(\nu(T'_2)-2)\leq 1$.
 Then $T_1+T_2$ is a Salem tree. (If $(\nu(T'_1)-2)(\nu(T'_2)-2)> 1$
then the reciprocal polynomial of $T_1+T_2$ has two roots outside
the unit circle.)
\item[(c)] Every Salem tree is of  type (a) or type (b).

\end{enumerate}
\end{theorem}

In case (a) of Theorem \ref{T-trees}, there is a single central
vertex  joined to $r$ cyclotomic subtrees $H_1$,
\ldots, $H_r$, while in case (b)  we have a central edge with each
endvertex joined to one or more cyclotomic subtrees
$H_1$, \ldots, $H_r$, $K_1$, \ldots, $K_s$:

\begin{picture}(200,200)(-100,-100)
\put(0,0){\circle*{5}} \put(0,0){\line(3,1){47.4}}
\put(0,0){\line(-3,1){47.4}} \put(0,0){\line(0,1){50}}
\put(45,-22){\circle*{2}} \put(40,-30){\circle*{2}}
\put(30,-35){\circle*{2}} \put(15,-39){\circle*{2}}
\put(0,-40){\circle*{2}} \put(-15,-39){\circle*{2}}
\put(-30,-35){\circle*{2}} \put(-40,-30){\circle*{2}}
\put(-45,-22){\circle*{2}} \put(0,65){\circle{30}}
\put(61.7,20.6){\circle{30}} \put(-61.7,20,6){\circle{30}}
%\put(-2,-10){$v$}
\put(-5,61){$H_1$} \put(56.7,16.6){$H_2$}
\put(-66.7,16.6){$H_r$}
\end{picture}

\begin{picture}(0,0)(-300,-120)
\put(0,0){\circle*{5}} \put(50,0){\circle*{5}}
\put(0,0){\line(1,0){50}} \put(0,0){\line(-1,3){15.8}}
\put(0,0){\line(-1,-3){15.8}} \put(50,0){\line(1,3){15.8}}
\put(50,0){\line(1,-3){15.8}} \put(97,10){\circle*{2}}
\put(94,20){\circle*{2}} \put(98,0){\circle*{2}}
\put(97,-10){\circle*{2}} \put(94,-20){\circle*{2}}
\put(-44,20){\circle*{2}} \put(-47,10){\circle*{2}}
\put(-48,0){\circle*{2}} \put(-47,-10){\circle*{2}}
\put(-44,-20){\circle*{2}} \put(-20.6,61.7){\circle{30}}
\put(-20.6,-61.7){\circle{30}} \put(70.6,61.7){\circle{30}}
\put(70.6,-61.7){\circle{30}}
%\put(2,-10){$v$} \put(42,-10){$w$}
\put(-25.6,57.7){$H_1$} \put(-25.6,-65.7){$H_r$}
\put(65.6,57.7){$K_1$} \put(65.6,-65.7){$K_s$}
\end{picture}

\begin{proof}
\begin{enumerate}
\item[(a)]
Take $\epsilon>0$ such that $R_T$ does not vanish on the interval
$I=(1,1+\epsilon)$.
Since $T'$ is cyclotomic, $R_{T'}>0$ on $(1,\infty)$, and hence
in particular $R_{T'}>0$ on $I$.
Since $\nu(T')>2$, $q_{T}(1)=1/(2-\nu(T'))<0$,
so $R_{T'}/R_T<0$ on $I$.
Hence $R_T<0$ on $I$.
Since $R_T(z)\to\infty$ as $z\to\infty$, $R_T$ has at least one
root on $(1,\infty)$.
By interlacing, $R_T$ cannot have more than one root on $(1,\infty)$,
since $R_{T'}$ has none.
This gives the first result.

\item[(b)] Let $T=T_1+T_2$.
Take $\epsilon>0$ such that neither $R_T$ nor $R_{T'}$ vanish
on $I=(1,1+\epsilon)$.
Now $T'$ is the forest $\{T_1',T_2\}$,
so that $R_{T'}$ has one root on $(1,\infty)$, a root of $R_{T_2}$.
By interlacing, $R_T$ has one or two roots on that interval.

On $I$, $R_{T'}<0$, and $z+1-zq_{T_2'}<0$, since $T_2$ is a Salem
tree of type (a).
If $(\nu(T_1')-2)(\nu(T_2')-2)<1$, then $R_{T'}/R_T>0$ on $I$,
so $R_T<0$ on $I$.
Since $R_T(z)\to\infty$ as $z\to\infty$, $R_T$ has an odd number
of roots, hence exactly one, on $(1,\infty)$.
(On the other hand, if $(\nu(T_1')-2)(\nu(T_2')-2)>1$, then
$R_{T'}/R_T<0$ on $I$, so $R_T>0$ on $I$, and then $R_T$ has an
even number, and hence two, roots on $(1,\infty)$.)

The only delicate case is if $(\nu(T'_1)-2)(\nu(T'_2)-2)- 1=0$.
Define the rational function
$f(z)=(z+1-z\pq_{T'_1})(z+1-z\pq_{T'_2})-z$, so that
$q_T=(z+1-zq_{T_2'})/f(z)$.
We need to identify the sign of $f(z)$ on $I$.
Putting $x=\sqrt{z}+1/\sqrt{z}$, the equation
$f(z)=0$ transforms to
$\psi(x):=(\sqrt{z}\pq_{T'_1}-x)(\sqrt{z}\pq_{T'_2}-x)=1$. Interlacing
implies that each
$\sqrt{z}\pq_{T'_i}$ (which is a function of $x$) is decreasing
between successive poles, and hence so too is each factor of
$\psi(x)$. But since $T_1$ and $T_2$ are Salem trees of type (a),
each factor of $\psi(x)$ is positive at
$x=2$; hence $\psi(x)<1$ as $x$ approaches 2 from above; hence
$f(z)<0$ on $I$.
Now, as before, we have $R_{T'}/R_T>0$ on $I$, so $R_T<0$ on $I$,
and the now familiar argument shows that $R_T$ has exactly one root
on $(1,\infty)$.

\item[(c)]

Suppose that $T$ is a tree such that $R_T$ has one root $>1$ but
is not of
 type (a).
Pick any vertex $t_0$ of $T$. Then, by interlacing, $T-\{t_0\}$ has
one component, $T_1$ say, that is a Salem tree, the other
components being cyclotomic. Let $t_1$ be the root of $T_1$ (the
vertex adjacent to $t_0$ in $T$). Now replace $t_0$ by $t_1$ and
repeat the argument, obtaining a new vertex $t_2$. If $t_2=t_0$ then
we are finished. Otherwise, we repeat the argument, obtaining a walk
on $T$, using vertices $t_0,t_1,t_2,\cdots$. Since $T$ has no
cycles, any walk in $T$  must eventually  double back on itself, so
that some $t_i$  equals $t_{i-2}$. Then $T$ is of the form
$T_1+T_2$, where $T_1$ and $T_2$ are of type (a), with roots
$t_{i-1}$ and $t_i$.
\end{enumerate}
\end{proof}

Note that while in case (a) $T$ is a rooted tree with the property
that removal of a single vertex gives a forest of cyclotomic
trees, in case (b) the tree $T_1+T_2$ has the property that
removal of the edge joining the roots of $T_1$ and $T_2$, with its
incident vertices, also gives a forest of cyclotomic trees.

Theorem \ref{T-trees} is a restriction of Theorem \ref{T-graphs}
(above) to trees. However, it is stronger, as we are able to say
precisely which trees are Salem trees. Theorem \ref{T-trees}  also
shows how to construct all Salem trees. To construct trees of type
(a), we take any collection of rooted cyclotomic trees, as listed in
Section \ref{S-quotients}, the sum of whose $\nu$-values exceeds
$2$. For trees of type (b), we take two such collections whose
$\nu$-values sum to $s_1$ and $s_2$ say, with $s_1\geq s_2>2$,
subject to the additional constraint that $(s_1-2)(s_2-2)\leq 1$. A
check on possible sums of $\nu$-values reveals that the smallest
possible value for $s_2>2$ is $85/42$, coming from the tree
$T(1,2,6)$, using the labelling of Figure \ref{F-TQ}. This implies
the upper bound $s_1\leq 44$ for $s_1$. Of course, when
$s_2>85/42$, the upper bound for $s_1$ will be smaller. Note too
that the condition $s_1\geq s_2$ implies that $s_2\leq 3$.

The first examples of Salem numbers of trace below $-1$ were
obtained using the construction in Theorem \ref{T-trees}(a) (see
\cite{ANTS}). The smallest known degree for a Salem number of trace
below $-1$ coming from a graph is of degree $460$, obtained when
$T'=\{A_{70}(1,69),D_{196}(182,14),D_{232}(220,12)\}$ in Theorem
\ref{T-trees}(a). Much smaller degrees have been obtained by other
means, and the minimal degree is known to be $20$ (\cite{ANTS}). It
is also now known that all integers occur as traces of Salem numbers
(\cite{MS}).

\subsection[]{ Earlier results} Theorem \ref{T-trees}(a) generalises
\cite[Corollary 9]{MRS}, which gave the same construction, but only for
starlike trees.
In 1988 Floyd and Plotnick \cite[Theorem 5.1]{FP}, without using graphs but
 using an unpublished result of Cannon, showed how to construct Salem numbers in a way
 equivalent to our construction using star-like trees. The same construction was also
 published by Cannon and Wagreich \cite[Proposition 5.2]{CW}
 and Parry \cite[Corollary 1.8]{P} in 1992. In 1999 Piroska Lakatos \cite[Theorem 1.2]{L1}
 deduced essentially the same star-like tree construction from results of A'Campo and Pena on Coxeter
 transformations. Also,  in 2001 Eriko Hironaka \cite[Proposition 2.1]{H}
  produced an equivalent construction,  in the context of knot theory,
 as the Alexander polynomial of a pretzel knot.

\section{Pisot graphs}\label{S-Pisot}

As we have seen in Section \ref{S-proof1}, a graph Pisot number is a
limit of graph Salem numbers whose graphs may be assumed to come
from a family obtained by taking a certain multigraph, and assuming
that some of its edges have an increasing number of subdivisions. We
use this family to define a graph having bi-coloured vertices: we
start with the multigraph, with black vertices. For every
increasingly subdivided pendant edge, we change the colour of the
pendant vertex to white, while for an increasing internal edge we
subdivide it with two white vertices. Thus a single white vertex
represents a pendant-increasing edge, while a pair of adjacent white
vertices represents an increasing internal edge. These Pisot graphs
in fact represent a sequence of Salem numbers tending to the
 Pisot number. Now, we have seen in the proof of Theorem \ref{T-closed} that the
  limit point of the Salem numbers corresponding to a Salem graph with an
 increasing internal edge is the same as that of the graph when this edge is
  broken in the middle. Hence for any Pisot graph we can remove any edge joining
   two white vertices without changing the corresponding Pisot number. (Doing
    this may disconnect the graph, in which case only one of the connected
     components corresponds to the Pisot number.) It follows
      that every graph Pisot number has a graph all of whose white vertices are pendant (have
      degree $1$).

For Pisot graphs that are trees ({\em Pisot trees}), and furthermore
have all white vertices pendant, we can define their quotients by
direct extension of the quotient of an ordinary tree (that is, one
without white vertices, as in Section  \ref{S-trees}). Now from
Section \ref{S-quotients} the path $A_n(1,n)$ has quotient
$(z^n-1)/(z^{n+1}-1)$, which, for $z>1$ tends to $1/z$ as
$n\to\infty$. Thus, following \cite[p. 315]{MRS}, we can take the
quotient of a white vertex $\circ$ to be $1/z$, and then calculate
the quotient of these trees in the same way as for ordinary trees.
The irreducible factor of its denominator with a root in $|z|>1$
then gives the minimal polynomial of the Pisot number.

\begin{figure}[h]
\begin{center}
\leavevmode
\hbox{%
\epsfxsize=3.3in
\epsffile{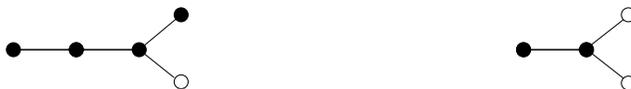}
}
\end{center}
\caption{  Pisot graphs of the smallest Pisot number (minimal
polynomial $z^3-z-1$), and for the smallest limit point of Pisot
numbers (minimal polynomial $z^2-z-1$). See end of Section
\ref{S-Pisot}.}\label{F-smallpisot}
\end{figure}

For instance, for the two Pisot trees in Figure \ref{F-smallpisot}, take their roots to be the central vertex.
Then we can use Lemma \ref{L-quotient}(i) to compute the quotient of the left-hand one to be
$$
\frac{1}{z+1-z\left(\frac1z+\frac{z-1}{z^2-1}+\frac{z^2-1}{z^3-1}\right)}=\frac{(z+1)(z^2+z+1)}{z(z^3-z-1)},
$$
so that the corresponding Pisot number has minimal polynomial $z^3-z-1$.
Similarly, the right-hand one has quotient $\frac{z+1}{z^2-z-1}$, with minimal polynomial $z^{2}-z-1$.

\section{Small elements of the derived sets of Pisot
numbers}\label{S-Bertin}

In this section we give a proof  of a graphical version of the following
result of Bertin \cite{Be}. Recall that the ($1$st) derived set of
a given real set is the set of limit points of the set, while for
$k\geq 2$ its $k$-th derived set is the set of limit points of its
$(k-1)$-th derived set.

\begin{theorem}\label{T-Bertin} Let $k\in\N$. Then $(k+\sqrt{k^2+4})/2$ belongs to the $(2k-1)$-th
derived set of the set $S_{\gra}$ of graph Pisot numbers, while
$k+1$ belongs to the $(2k)$-th derived set of $S_{\gra}$.
\end{theorem}

Bertin's result was that these numbers belonged to the corresponding derived set of $S$, rather than that of
$S_{\gra}$. They are the smallest known elements of the relevant derived set of $S$.

\begin{figure}[h]
\begin{center}
\leavevmode
\hbox{%
\epsfxsize=2.1in
\epsffile{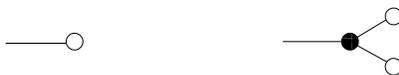}
}
\end{center}
\caption{The subtrees used to make small elements of the derived
sets of the set of graph Pisot numbers. Their Pisot quotients are
$1/z$ and $1/(z-1)$ respectively. They give such elements as a
limit of increasing graph Pisot numbers. See Theorem
\ref{T-Bertin}. }\label{F-Bertinbits}
\end{figure}

\begin{proof}
The proof consists simply of exhibiting two families of trees
containing $2k$ and $2k+1$  white vertices respectively, and showing
that their reciprocal polynomials are $z^2-kz-1$ and $z-(k+1)$. From
the discussion above, this will show that their zeros in $|z|>1$,
namely those given in the statement of the Theorem, are in the
$(2k-1)$-th and $(2k)$-th derived set  of the set of Pisot numbers,
respectively. For the graph with $2k$  white vertices we take $k$ of
the $3$-vertex graphs shown in Figure \ref{F-Bertinbits} joined to a
central vertex, while for the graph with $2k+1$ vertices we take the
same graph with  one extra white vertex joined to the central vertex
(the other graph shown in this figure). The result is shown in
Figure \ref{F-Bertin} for $k=5$. We can use Lemma \ref{L-quotient},
extended to include trees containing an infinite path. This shows
that the tree has quotient $(z+1-kz/(z-1))^{-1}=(z-1)/(z^2-kz-1)$
when it has $2k$ white vertices, and quotient
$(z+1-z(k/(z-1)+1/z))^{-1}=(z-1)/(z(z-(k+1)))$ when it has $2k+1$
white vertices. The poles of these quotients give the required Pisot
numbers.

\end{proof}

\begin{figure}[h]
\begin{center}
\leavevmode
\hbox{%
\epsfxsize=3.3in
\epsffile{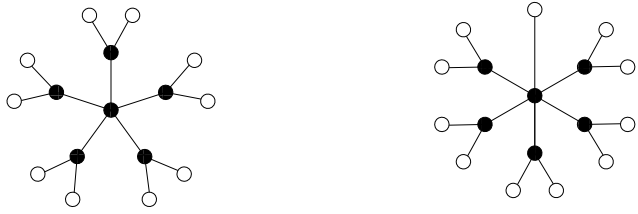}
}
\end{center}
\caption{The infinite graphs showing that $(k+\sqrt{k^2+4})/2$ belongs to the $(2k-1)$-th
derived set of the set of graph Pisot numbers (left, $k=5$ shown), and that  $k+1$ belongs to the
 $(2k)$-th
derived set (right). Here increasing sequences are produced---see Theorem \ref{T-Bertin}.
}\label{F-Bertin}
\end{figure}

The graphs of Figure \ref{F-Bertin} show how the elements of the
derived sets are limits from below of elements of $S_{\gra}$. We can
also show that they are limits from above, using the $5$- and
$11$-vertex graphs of Figure \ref{F-Bertinbits_above} to construct
Pisot graphs showing these numbers to be elements of the relevant
derived set by showing them to be limit points from above rather
than below. The graphs in Figure \ref{F-Bertin_above} are examples
of this construction. Further, one could construct graphs using a
mixture of subgraphs from Figures \ref{F-Bertinbits} and
\ref{F-Bertinbits_above}. Thus, if we distinguished two types of
limit point depending on whether the point was a limit from below or
from above, we could define two types of derived set, and hence, by
iteration, an $(n_-,n_+)$-derived set of $S_{\gra}$. This mixed
construction would produce elements of these sets.

\begin{figure}[h]
\begin{center}
\leavevmode
\hbox{%
\epsfxsize=3.1in
\epsffile{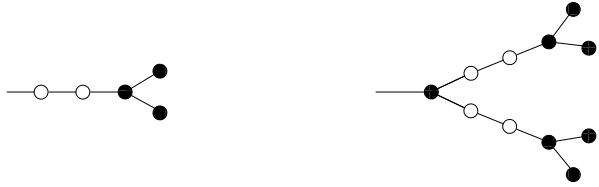}
}
\end{center}
\caption{The subtrees used to make small elements of the derived sets of the set of
graph Pisot numbers. Their Pisot quotients are $1/z$ (left) and $1/(z-1)$ (right).
They give such elements as a limit of decreasing graph Pisot numbers. See the remarks after Theorem \ref{T-Bertin}.
}\label{F-Bertinbits_above}
\end{figure}

\begin{figure}[h]
\begin{center}
\leavevmode
\hbox{%
\epsfxsize=5.3in
\epsffile{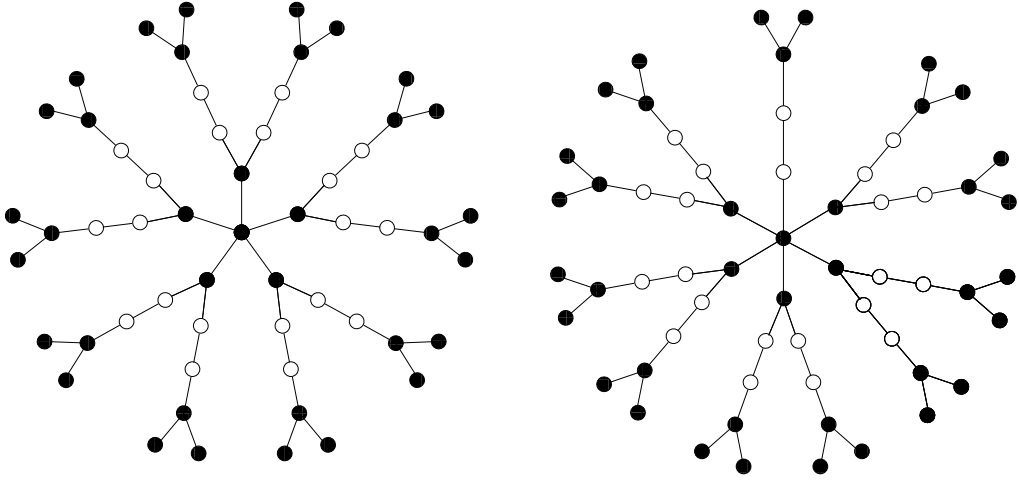}
}
\end{center}
\caption{The infinite graphs showing that $(k+\sqrt{k^2+4})/2$ belongs to the $(2k-1)$-th
derived set of the set of graph Pisot numbers (left, $k=5$ shown), and that  $k+1$ belongs to the
 $(2k)$-th
derived set (right). Here decreasing sequences are produced---see the remarks after Theorem \ref{T-Bertin}.
}\label{F-Bertin_above}
\end{figure}

\section{The Mahler measure of graphs}\label{S-Mahler}

In this section we find (Theorem \ref{T-strong}) all the graphs of Mahler measure less than
$\rho:=\frac{1}{2}(1+\sqrt{5})$.  Our
definition of Mahler measure for graphs---see below---seems natural. This is because we then obtain
as a corollary that the strong version of
``Lehmer's Conjecture'', which states that $\tau_1$ is the smallest Mahler measure greater than $1$
of any algebraic number,
is true for graphs:

\begin{corollary} The Mahler measure of a graph is either $1$ or at least $\tau_1=1.176280818\cdots$,
the largest real zero of
Lehmer's polynomial $L(z)=z^{10}+z^9-z^7-z^6-z^5-z^4-z^3+z+1$. Among
connected graphs, this minimum Mahler measure is attained only for
the graph $T(1,2,6)$ defined in Figure \ref{F-TQ} ($=$ the
Coxeter graph $E_{10})$.
\end{corollary}

 We define the {\em Mahler measure} $M(G)$ of an $n$-vertex graph $G$ to
be $M(z^n\pchi_G(z+1/z))$, where $\pchi_G$ is the characteristic
polynomial of its adjacency matrix, and $M$ of a polynomial also
denotes its Mahler measure.  Recall that for a monic polynomial
$P(z)=\prod_i(z-\alpha_i)$ its Mahler measure is defined to be
$M(P)=\prod_i \max(1,|\alpha_i|)$. When $G$ is bipartite, $M(G)$ is
also the Mahler measure of its reciprocal polynomial
$R_G(z)=z^{n/2}\pchi_G(\sqrt{z}+1/\sqrt{z})$. This is because then
$M(z^n\pchi_G(z+1/z))=M(R_G(z^2))$.

The graphs having Mahler measure $1$ are precisely the cyclotomic graphs.

It turns out that  the connected graphs of smallest Mahler measure
bigger than $1$ are all trees.
Using the notation of \cite{CR}, define the trees $T(a,b,c)$ and
$Q(a,b,c)$ as in Figure \ref{F-TQ}.

\begin{theorem} \label{T-strong}If $G$ is a connected graph whose Mahler measure $M(G)$ lies in the interval
$(1,\rho )$ then $G$ is one of the following trees:
\begin{itemize}
\item $G=T(a,b,c)$ for $a\leq b\leq c$ and
    \begin{itemize}
    \item[] $a=1,b=2,c\geq 6$
    \item[] $a=1,b\geq 3,c\geq 4$
    \item[] $a=2,b=2,c\geq 3$
    \item[] $a=2,b=3,c=3$
    \end{itemize}
or
\item $G=Q(a,b,c)$ for $a\leq c$ and
    \begin{itemize}
    \item[] $a=2,b\geq 1,c=3$
    \item[] $a=2,b\geq 3, 4\leq c\leq b+1$
    \item[] $a=3,4\leq b\leq 13,c=3$
    \item[] $a=3,5\leq b\leq 10,c=4$
    \item[] $a=3,7\leq b\leq 9,c=5$
    \item[] $a=3,8\leq b\leq 9,c=6$
    \item[] $a=4,7\leq b\leq 8,c=4$.
    \end{itemize}
\end{itemize}
All these graphs $G$ are all Salem graphs, apart from $Q(3,13,3)$,
whose polynomial has two zeros on $(1,2)$, so that all $M(G)$ apart
from $M(Q(3,13,3))$ are Salem numbers.
 Also, the set of limit points of the
set of all $M(G)$ in $(1,\rho\, ]$ consists of the  graph  Pisot number $\rho $ and the
  graph  Pisot numbers that are zeros of $z^k(z^2-z-1)+1$ for $k=2,3,\dots$, which approach $\rho $ as $k\to\infty$.

Furthermore, the only $M(G)<1.3$ are $M(T(1,2,c))$ for $c=6,7,8,9,10$, these values increasing with $c$. (Also
$M(T(1,2,9))= M(T(1,3,4))$.)
\end{theorem}

\begin{figure}[h]
\begin{center}
\leavevmode
\psfragscanon
\psfrag{a}[l]{$\overbrace{\phantom{abcdefghijk}}^a$}
\psfrag{b}[l]{$\overbrace{\phantom{abcdefghijk}}^b$}
\psfrag{c}[l]{$\overbrace{\phantom{abcdefghijk}}^c$}
\psfrag{T}{$T(a,b,c)$}
\psfrag{Q}{$Q(a,b,c)$}
\psfrag{a1}[l]{$\overbrace{\phantom{abcdefghijk}}^{a-1}$}
\psfrag{b1}[l]{$\overbrace{\phantom{abcdefghijk}}^{b-1}$}
\psfrag{c1}[l]{$\overbrace{\phantom{abcdefghijk}}^{c-1}$}
\includegraphics[scale=0.6]{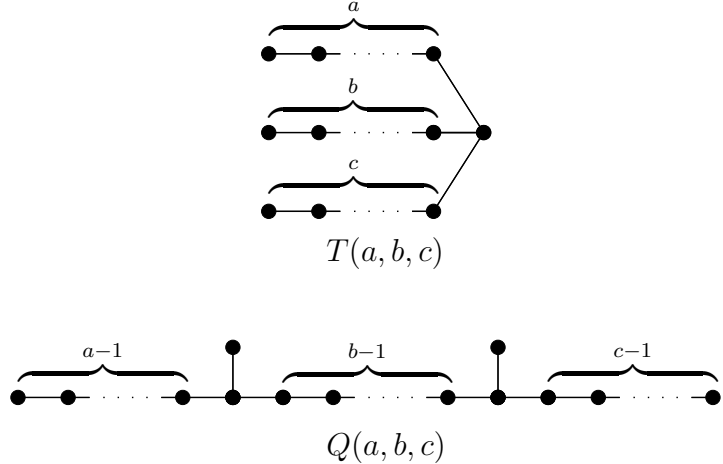}
\end{center}
\caption{The trees $T(a,b,c)$ and $Q(a,b,c)$ }\label{F-TQ}
\end{figure}

In \cite[Theorem 1.1]{H}, Hironaka shows essentially that Lehmer's number
 $\tau_1=M(T(1,2,6))$ is the smallest Mahler measure of a starlike tree.
 The graph Pisot numbers in the Theorem have been shown by Hoffman (\cite{H2}) in 1972 to
 be limit points of transformed graph indices, which is equivalent to our
 representation of them as limits of graph Salem numbers.

It is clear how to extend the theorem to nonconnected graphs:
since $\tau_1^3>\rho$, for such a
graph $G$ to have $M(G)\in(1,\rho )$,  one or two connected
components must be as described in the theorem, with all other
connected components cyclotomic. Using the results of the
Theorem, it is an easy exercise to check the possibilities.

\begin{proof} The proof depends heavily on results of Brouwer and Neumaier \cite{BN} and Cvetkovi\' c, Doob and Gutman
\cite{CDG}, as described conveniently by Cvetkovi\'c and Rowlinson
in their survey paper \cite[Theorem 2.4]{CR}. These results tell
us precisely whch connected graphs have largest eigenvalue in the
interval $(2,\sqrt{2+\sqrt{5}}\,]=(2,2.058\cdots]$. They are all
trees of the form $T(a,b,c)$ or $Q(a,b,c)$.  Those of the form
$T(a,b,c)$ are precisely those given in the statement of the
theorem. As they are starlike trees, they have, by \cite[Lemma
8]{MRS}, one eigenvalue $\la>2$, and so their reciprocal
polynomial $R_{T(a,b,c)}$ has a single zero $\be$ on $(1,\infty)$
with $\be^{1/2}+\be^{-1/2}=\la$, and $M(T(a,b,c))=\be$. Since
$\rho^{1/2}+\rho^{-1/2}=\sqrt{2+\sqrt{5}}$, we have $\be\in(1,\rho
]$.

 From the previous paragraph it is clear that
\begin{itemize}
\item all graphs $G$ with exactly one eigenvalue in $(2,\sqrt{2+\sqrt{5}}\,]$  have $M(G)\in(1,\rho )$;
\item only graphs $G$ with largest eigenvalue in $(2,\sqrt{2+\sqrt{5}}\,]$ can have $M(G)\in(1,\rho )$.
\end{itemize}

 It remains to see which of the graphs $Q(a,b,c)$ having largest eigenvalue in this interval
actually do have $M(G)\in(1,\rho\, ]$. The graphs of this type
given in the Theorem are all those with one eigenvalue in
$(2,\sqrt{2+\sqrt{5}}\,]$,  along with $Q(3,13,3)$ which, although
having two eigenvalues greater than $2$, nevertheless does have
$M(G)<\rho $. The other graphs with largest eigenvalue in
$(2,\sqrt{2+\sqrt{5}}\,]$ are, from the theorem cited above:

$Q(3,b,3)$ for $b\geq 14$,

$Q(3,b,4)$ for $b\geq 11$,

$Q(3,b,5)$ for $b\geq 10$,

$Q(3,b,6)$ for $b\geq 10$,

$Q(3,b,c)$ for $b\geq c+2,c\geq 7$,

$Q(4,b,4)$ for $b\geq 9$,

$Q(4,b,5)$ for $b\geq 8$,

$Q(4,b,c)$ for $b\geq c+4,c\geq 6$,

$Q(a,b,c)$ for $a\geq 5,b\geq a+c,c\geq 5$.

We must show that none of these trees $G$ have $M(G)\leq\rho$.
We can reduce this infinite list  to a small finite one by the
following simple observation. Suppose we remove the $k$-th vertex
from the central path of the tree $Q(a,b,c)$, splitting it into
$T(1,a-1,k-1)$ and $T(1,c-1,b-1-k)$. By interlacing we have, for
$k=2,\dots,b-2$,
\begin{equation}\label{E-two-T}
M(Q(a,b,c))\geq M(T(1,a-1,k-1))\times M(T(1,c-1,b-1-k)).
\end{equation}
Now
\begin{align*}
M(T(1,2,6))=&1.176280818\cdots,\\
 M(T(1,2,9))=& M(T(1,3,4))= 1.280638156\cdots,\\
  M(T(1,3,6))=& M(T(1,4,4))= 1.401268368\cdots,
  \end{align*}
   from
which we have that both of $M(T(1,2,6))\times M(T(1,3,6))=M(T(1,2,6))\times M(T(1,4,4))$ and
$M(T(1,2,9))^2=M(T(1,3,4))^2$ are greater than $\rho $. Since $M(T(a,b,c))$ is, when greater than $1$,
an increasing function of $a$, $b$ and $c$ separately, and, of course, independent of the order of $a,b,c$, we can show
that all but $18$ of the above $Q(a,b,c)$ have $M(Q(a,b,c))>\rho $. Applying (\ref{E-two-T}), we have
\begin{itemize}
\item $M(Q(3,b,3))\geq M(T(1,2,9))\times M(T(1,2,b-11))>\rho $ for $b\geq 20$. Cases $b=14\cc 19$ must be
checked individually.
\item $M(Q(3,b,4))\geq M(T(1,2,6))\times M(T(1,3,b-8))>\rho $ for $b\geq 14$. Check $b=11,12,13$ individually.
\item $M(Q(3,b,5))\geq M(T(1,2,6))\times M(T(1,4,b-8))>\rho $ for $b\geq 12$. Check $b=10,11$.
\item $M(Q(3,b,6))\geq M(T(1,2,6))\times M(T(1,5,b-8))> M(T(1,2,6))\times M(T(1,4,b-8))>\rho $ for $b\geq 12$.
 Check $b=10,11$.
\item $M(Q(3,b,7))\geq M(T(1,2,6))\times M(T(1,6,b-8))>\rho $ for $b\geq 11$.
 Check $b=9,10$.
\item $M(Q(3,b,8))\geq M(T(1,2,6))\times M(T(1,7,b-8))> M(T(1,2,6))\times M(T(1,6,b-8))>\rho $ for $b\geq 11$.
 Check $b=10$.
\item For $c\geq 9$, $M(Q(3,b,c))\geq M(T(1,2,6))\times M(T(1,c-1,b-8))>\rho $ for $b\geq 11$.
\item $M(Q(4,b,4))\geq M(T(1,3,4))\times M(T(1,3,b-6))>\rho $ for $b\geq 10$.
Check $b=9$.
\item $M(Q(4,b,5))\geq M(T(1,3,4))\times M(T(1,4,b-6))>\rho $ for $b\geq 9$.
Check $b=8$.
\item For $c\geq 6$, $M(Q(4,b,c))\geq M(T(1,3,4))\times M(T(1,c-1,b-6))>M(T(1,3,4))\times M(T(1,4,b-6))
>\rho $ for
$b\geq 9$.
\item For $a\geq 5$, $c\geq 5$, $M(Q(a,b,c))\geq M(T(1,a-1,3))\times M(T(1,c-1,b-5))\geq M(T(1,4,3))
\times M(T(1,4,b-5))>\rho $ for
$b\geq 8$.
\end{itemize}

We remark that it is straightforward, with computer assistance, using Lemma \ref{L-quotient}, to make the checks required
 in the proof. Denoting
by $q_k(a,b,c)$ the quotient of $Q(a,b,c)$ with root at the $k$-th vertex of the central path,
 and by $t(a,b,c)$ the
quotient of $T(a,b,c)$ having root at the endvertex of the
$c$-path, this lemma tells us that
\begin{align*}
q_k(a,b,c)= &(z+1-z(t(1,a-1,k-1)+t(1,c-1,b-1-k)))^{-1},\\
t(a,b,c)=& (z+1-zt(a,b,c-1))^{-1}
\end{align*}
with $t(a,b,0)=\frac{(z^{a+1}-1)(z^{b+1}-1)}{(z-1)(z^{a+b+2}-1)}$,
the quotient of the rooted path $A_{a+b+1}(a+1,b+1).$ Then the
denominator of $q_{k}(a,b,c)$ gives the reciprocal polynomial of
$Q(a,b,c)$, at least up to a cyclotomic factor (one can show using Lemma 
  \ref{L-addedge}(i) and Lemma \ref{L-quotient}(i) that all roots $>1$ of the reciprocal polynomial of $Q(a,b,c)$ are indeed poles of its quotient).
  
   Concerning the limit points of $M(G)\cap [1,\rho)$, one
can check that
\begin{itemize}
\item $M(T(1,b,c))\to M(z^b(z^2-z-1)+1)$ as $c\to \infty$;

\item $M(T(2,2,c))\to \rho$ as $c\to \infty$;
\item For $c\geq 3$, $ M(Q(2,b,c))\to M(z^{c-1}(z^2-z-1)+1)$ as $b\to
\infty$.
\end{itemize}
Of course, by Salem's classical construction,
$M(z^b(z^2-z-1)+1)\to \rho$ as $b\to\infty$. Note too that
$M(z^{2}(z^2-z-1)+1)=M(z^3-z-1)$, the smallest Pisot number.

\end{proof}

From the proof, and the fact that all Pisot numbers in $[1,\rho)$
are known
 (see Bertin {\em et al} \cite[p.133]{BDGPS}) we have the following.

\begin{corollary}\label{C-small-Pisot} The only graph Pisot numbers in $[1,\rho)$ are the roots of
$z^n(z^2-z-1)+1$ for $n\geq 2$. The other Pisot numbers in this interval,
namely the roots of $z^6-2z^5+z^4-z^2+z-1$ and of
 $z^n(z^2-z-1)+z^2-1$ for $n\geq 2$, are not graph Pisot numbers.
\end{corollary}

\section{Small Salem numbers from graphs}\label{S-SmSa}
The notation $\tau_n$ indicates the $n$th Salem number
 in Mossinghoff's table \cite{M}, listing all 47 known
Salem numbers that are smaller than $1.3$. (This is an update of that in \cite{Bo}.)

We have seen that the only numbers in this list that are elements of
$T_{\gra}$ are $\tau_1$, $\tau_7$, $\tau_{19}$, $\tau_{23}$, and
$\tau_{41}$. On the other hand, if we apply the construction in
Theorem \ref{T-trees}(b) with $T_1=T_2$, then from the explicit
formula in Lemma \ref{L-quotient}(ii) we see that the Salem number
produced is automatically the square of a smaller Salem number. ( If
$\pq_{T_1}=q/p$ in lowest terms, then from Lemma
\ref{L-quotient}(ii) we have that the squarefree part of
$R_{T_1+T_1}$ is $f(z)=q(z)^2-zp(z)^2$. Now
$f(z^2)=(q(z^2)-zp(z^2))(q(z^2)+zp(z^2))$, and the gcd of these two
factors divides $z$.  Hence, if $\tau\ne0$ and $f(\tau)=0$, then
$\sqrt{\tau}$ and $-\sqrt{\tau}$ are roots of the different factors
of $f(z^2)$ and are not algebraic conjugates. In particular, if
$\tau$ is a Salem number then so is $\sqrt{\tau}$.  We can apply
this construction regardless of the value of the quotient of $T_1$.)
In this way we can produce $\tau_2^2$, $\tau_3^2$, $\tau_5^2$,
$\tau_{12}^2$, $\tau_{21}^2$, $\tau_{23}^2$, and $\tau_{41}^2$ as
elements of $T_{\gra}$.

These results, and a wider search for small powers of small
Salem numbers, are recorded in the following table.
A list of cyclotomic graphs indicates the components of
$T'$ in the construction of Theorem \ref{T-trees}(a);
two lists separated by a semi-colon indicate the components
of $T'_1$ and $T'_2$ in the construction of Theorem \ref{T-trees}(b).

\renewcommand{\arraystretch}{1.1}
\[
\begin{array}{|c|c|}
\hline
\rm Salem\ number & \rm Cyclotomic\  graphs \\ \hline
\tau_1 & D_9(0) \\
\tau_1^2 & D_{11}(3,8) \\
\tau_1^6 & E_7(1), \tilde{D}_4(0);A_5(2,4)\\
\tau_1^8 & E_6(1),A_2(1,2);E_7(5),\tilde{E}_6(3)\\ \hline
\tau_2^2 & E_8(7);E_8(7) \\ \hline
\tau_3^2 & E_7(6);E_7(6) \\
\tau_3^5 & A_1(1,1),A_9(2,8);D_{15}(8,7)\\ \hline
\tau_4^5 & E_6(4),D_7(1,6);D_{13}(3,10) \\ \hline
\tau_5^2 & E_6(1);E_6(1) \\
\tau_5^3 & E_6(1);\tilde{E}_8(7) \\
\tau_5^4 & E_6(4),D_{18}(12,6) \\
\tau_5^5 & A_4(1,4),A_4(1,4);D_4(1,3),D_8(1,7) \\
\tau_5^6 & A_1(1,1),A_3(2,2);D_6(2,4),D_8(4,4) \\ \hline
\tau_7 & D_{10}(0) \\
\tau_7^4 & E_6(1),A_1(1,1);E_6(1),A_1(1,1) \\
\tau_7^5 & A_7(2,6);D_4(1,3),\tilde{D}_{10}(5,5) \\
\tau_7^6 & E_7(3),D_7(4,3);D_9(1,8)\\ \hline
\tau_{10}^3   &   E_8(8);D_8(0)\\ \hline
\tau_{12}^2 & D_5(0);D_5(0) \\
\tau_{12}^3 & E_7(5);E_7(6) \\
\tau_{12}^5 & E_7(4),\tilde{E}_6(1);A_7(3,5)\\ \hline
\tau_{15}^2    &  D_{18}(6,12)\\
\tau_{15}^4 & A_1(1,1),D_{10}(0); A_1(1,1),D_{10}(0) \\ \hline
\tau_{16}^4 & E_7(1),D_9(1,8),D_8(2,6) \\ \hline
\tau_{19} & D_{11}(0) \\
\tau_{19}^3     & \tilde{E}_8(8);D_4(2,2)\\
\tau_{19}^4 & E_6(4),A_1(1,1);E_6(4),A_1(1,1)\\
\tau_{19}^5 & \tilde{E}_6(2),A_3(2,2);A_3(1,3),D_6(1,5)\\ \hline
\tau_{21}^2 & E_7(1);E_7(1) \\
\tau_{21}^5 & E_6(3),A_4(2,3);A_6(1,6) \\ \hline
\tau_{23} & E_8(1) \\
\tau_{23}^2 & \tilde{E}_8(6) \\
\tau_{23}^3     & E_7(2);D_6(1,5)\\
\tau_{23}^4 & \tilde{E}_7(3),D_{12}(9,3) \\ \hline
\tau_{35}^4 & E_6(4),E_7(1);A_2(1,2),A_6(1,6) \\ \hline
\tau_{41} & D_{13}(0) \\
\tau_{41}^2 & D_6(0);D_6(0) \\
\tau_{41}^3     & A_7(2,6);D_{10}(5,5)\\
\tau_{41}^4 & A_2(1,2),A_2(1,2);A_6(2,5),D_5(1,4)\\
\hline
\end{array}
\]

\end{document}